\documentclass[10pt,twoside,fleqn]{article}

\usepackage{amsmath, amsthm, amssymb, mathrsfs}
\usepackage{fancyhdr}
\usepackage{epsfig, graphicx}

\newtheorem{thm}{Theorem}[section]
\newtheorem{cor}[thm]{Corollary}

\newtheorem{lem}[thm]{Lemma}

\renewcommand{\theequation}{\arabic{section}.\arabic{equation}}

\newcommand{\T}    {\mathbb{T}}

\newcommand{\R}    {\mathbb{R}}
\newcommand{\C}    {\mathbb{C}}
\newcommand{\N}    {\mathbb{N}}
\newcommand{\Z}    {\mathbb{Z}}

\newcommand{\A}    {\mathbf{A}}
\newcommand{\CoM}  {\mathbf{BVT}}

\newcommand{\Pc}   {{\cal P}}

\newcommand{\K}    {{\cal K}}

\newcommand{\supp} {\textnormal{supp}}
\newcommand{\cp}   {\textnormal{cap}}
\newcommand{\dist} {\textnormal{dist}}
\newcommand{\diam} {\textnormal{diam}}
\newcommand{\Arg}  {\textnormal{Arg}}
\newcommand{\ang}  {\textnormal{Angle}}
\newcommand{\re}   {\textnormal{Re}}
\newcommand{\im}   {\textnormal{Im}}
\newcommand{\Log}  {\textnormal{Log}}

\newcommand{\spp}  {S}
\newcommand{\rspp} {\spp^\prime}
\newcommand{\wspp} {\widetilde\spp}
\newcommand{\mes}  {\lambda}

\newcommand{\Lm}   {\Lambda}
\newcommand{\arm}  {\varphi}
\newcommand{\angs} {\theta}

\newcommand{\cws}  {\stackrel{*}{\rightarrow}}
\newcommand{\cic}  {\stackrel{\scriptsize\cp}{\rightarrow}}

\newcommand{\boite}{\mbox{} \hfill \mbox{\rule{2mm}{2mm}}}

\begin{document}

\title{\large\bf Multipoint Pad\'e Approximants to Complex Cauchy Transforms with Polar Singularities}

\author{L. Baratchart and M. Yattselev}

\date{\normalsize \today}

\maketitle

\thispagestyle{empty}

\begin{abstract}
We study diagonal multipoint Pad\'e approximants to functions of the form
\[F(z) = \int\frac{d\mes(t)}{z-t}+R(z),\]
where $R$ is a rational function and $\mes$ is a complex measure with compact regular support included in $\R$, whose argument has bounded variation on the support. Assuming that interpolation sets are such that their normalized counting measures converge sufficiently fast in the  weak-star sense to some conjugate-symmetric distribution $\sigma$, we show that the counting measures of poles of the approximants converge to $\widehat\sigma$, the balayage of $\sigma$ onto the support of $\mes$, in the weak$^*$ sense, that the approximants themselves converge in capacity to $F$ outside the support of $\mes$, and that the poles of $R$ attract at least as many poles of the approximants as their multiplicity and not much more. 
\smallskip
\newline
{\em AMS Classification (MSC2000):} primary 41A20, 41A30, 42C05; secondary 30D50, 30D55, 30E10, 31A15.
\end{abstract}

\refstepcounter{section}
\section*{\normalsize\centering\thesection~ Introduction}
\label{sec:intro}

This paper is concerned with the asymptotic behavior of diagonal multipoint Pad\'e approximants to functions of the form
\begin{equation}
\label{fH}
F(z)=\int\frac{d\mes(t)}{z-t}+R(z),
\end{equation}
where $R$ is a rational function holomorphic at infinity and $\mes$ is a complex measure compactly and regularly supported on the real line. 

Diagonal multipoint Pad\'e approximants are rational interpolants of type $(n,n)$ where, for each $n$, a set of $2n+1$ interpolation points has been prescribed, one of which is infinity. Moreover, we assume that the interpolation points converge sufficiently fast to a conjugate-symmetric limit distribution whose support is disjoint from both the poles of $R$ and the convex hull of $\supp(\lambda)$, the support of $\mes$ (see (\ref{eq:suffast})).

To put our results into
perspective, let us begin with an account of the existing literature.
When $\mes$ is a positive
measure and $R\equiv0$ (in this case $F$ is referred to as a Markov function),
the study of diagonal Pad\'e approximants to $F$ at infinity 
goes back to A. A. Markov who showed
(see \cite{Mar95}) that they converge uniformly to 
$F$ on compact subsets of $\overline\C\setminus I$, where $I$ is
the convex hull of $\supp(\mes)$. 
Later this work was extended to multipoint Pad\'e approximants with
conjugate-symmetric interpolation schemes
by A. A. Gonchar and G. L\'opez Lagomasino 
in \cite{GL78}. A cornerstone of the theory is the close relationship between 
 Pad\'e  approximants to Markov functions and orthogonal polynomials,
since the denominator of the $n$-th diagonal approximant is the $n$-th
orthogonal polynomial in $L^2(d\mes)$ (resp. $L^2(d\mes/p)$, where $p$ is a
polynomial vanishing at finite interpolation points). For further references
and sharp error rates, we refer the reader to the monographs
\cite{StahlTotik,Totik}.

Another generalization of Markov's result was obtained
by A. A. Gonchar on adding polar singularities, i.e. on making
$R\not\equiv0$. He proved in \cite{Gon75a} that Pad\'e approximants still converge to 
$F$ locally uniformly in $\overline\C\setminus(\rspp\cup I)$, where $\rspp$ is 
the set of poles of $R$, provided 
that $\mes$ is a positive measure with singular part supported on a set of 
logarithmic capacity zero. Subsequently, it was shown by E. A. Rakhmanov 
in \cite{Rakh77b} 
that weaker assumptions on $\mes$ can spoil the convergence, but at 
the same time that if the coefficients of $R$ are real, then the locally
uniform convergence holds for any positive  $\mes$. Although it is
not a concern to us here,
let us mention that one may also relax
the assumption that $\supp(\mes)$ is compact. In particular,
Pad\'e and multipoint Pad\'e approximants to
Cauchy transforms of positive measures supported in $[0,\infty]$ 
(such functions are said to be of Stieltjes type) were investigated by 
G. L\'opez Lagomasino in \cite{LL78,LL88}. Let us also stress that polynomials satisfying certain Sobolev-type orthogonality exhibit an asymptotic behavior quite similar to that of the denominators of diagonal Pad\'e approximants to functions of the form (\ref{fH}) with non-trivial $R$ \cite{LLMarVA95}.

The case of a complex measure was taken up by G. Baxter in \cite{Bax61}
and by J. Nuttall and S. R. Singh in \cite{NS80}, who established strong
asymptotics of non-Hermitian orthogonal polynomials
on a segment for measures that are absolutely continuous with respect to the
(logarithmic) equilibrium distribution of that segment, 
and whose density satisfy appropriate
conditions expressing, in one way or another, that it is smoothly invertible. 
These results entail
that the Pad\'e approximants to $F$ converge uniformly to the latter 
on compact subsets of $\overline\C\setminus I$ when 
$R\equiv0$ and $d\mes/d\mu_I$ meets these
conditions (here $\mu_I$ indicates the equilibrium distribution on $I$).
For instance Baxter's condition is that 
$\log d\mes/d\mu_I$, when extended periodically, has an absolutely 
summable Fourier series. When $d\mes(t)/dt$ is holomorphic and nonvanishing
on a neighborhood of $I$, still stronger asymptotics, which apply to 
multipoint Pad\'e approximants as well, were recently obtained by
A. I. Aptekarev in \cite{Ap02} (see also \cite{AVA04}), using the matrix Riemann-Hilbert
approach pioneered by P. Deift and X. Zhou (see e.g. \cite{Deift}). Eventhough it is not directly related to the present work we mention for completeness another approach to analyzing the asymptotics of Pad\'e approximants based on three term recurrence relations \cite{BarrLLTorr95}.

Meanwhile H. Stahl opened up new perspectives in his
pathbreaking papers \cite{St89,St97}, where he studied diagonal
Pad\'e approximants
to (branches of) multiple-valued functions that can be continued 
analytically without restriction except over a set of capacity zero 
(typical examples are functions
with poles and branchpoints). By essentially representing the ``main''
singular part of the function as a Cauchy
integral over a system of cuts of minimal capacity,
and through a deep analysis of zeros of non-Hermitian 
orthogonal polynomials on such systems of cuts, he established the 
asymptotic distribution of
poles and subsequently the convergence in capacity of the Pad\'e 
approximants on the
complement of the cuts. In \cite{GRakh87} this construction was 
generalized to certain carefully chosen 
multipoint Pad\'e approximants by A. A. Gonchar and E. A. Rakhmanov, 
who, in particular, used it to illustrate the sharpness of O. G. Parfenov's theorem (formerly Gonchar's 
conjecture)
on the rate of approximation by rational functions over compact subsets of the
domain of holomorphy, see \cite{Par86}. Of course the true power of this method lies with the fact that it allows one
to deal with measures supported on more general systems of arcs than a
segment, which is beyond the scope of the present paper. However, since a 
segment is the simplest example of an arc of minimal logarithmic
capacity connecting
two points, the results we just mentioned apply in
particular to functions of the form (\ref{fH}), where $\mes$ is a complex
measure supported
on a segment which is absolutely continuous there with continuous density 
that does not vanish outside a set of capacity zero. 
By different, operator-theoretic methods, combined with a well-known
theorem of E. A. Rakhmanov on ratio asymptotics (see \cite{Rakh77a}),
A. Magnus further showed
that the diagonal Pad\'e approximants to $F$ converge \emph{uniformly} 
on compact subsets of $\overline\C\setminus I$
when $R\equiv0$ and $d\mes/dt$ is non-zero almost everywhere with 
continuous argument \cite{Mag87}. The existence of a uniformly convergent subsequence of diagonal Pad\'e approximants to (\ref{fH}) with non-trivial $R$ was shown in \cite{Suet00} whenever $\supp(\mes)$ is a disjoint union of analytic arcs in ``general position'' of minimal capacity and $d\mes/dt$ is sufficiently smooth and non-vanishing. Moreover, when $\supp(\mes)$ is a union of several intervals and the density of the measure is real analytic, the behavior of the zeros that do not approach $\supp(\lambda)$ nor the poles of $R$ can be described by the generalized Dubrovin system of non-linear differential equations \cite{Suet02}.

In contrast with previous work, the present approach allows the complex measure $\mes$ to vanish on a large subset of $I$. Specifically, we require that the total variation measure $|\mes|$ has compact regular support and that it is not too thin, say, larger than a power of the radius on relative balls \emph{of the support} (see the definition of the class  $\CoM$ in Section \ref{sec:pade}). In particular, this entails that $\supp(\mes)$ could be a thick Cantor set, or else the closure of a union of infinitely many intervals; such cases could not be handled by previously known methods. Although fairly general, these conditions could be further weakened, for instance down to the
$\Lambda$-criterion introduced by H. Stahl and V. Totik 
in \cite{StahlTotik}\footnote{This depends on the corresponding generalization of 
the results in \cite{BKT05} to be found in \cite{thKus}, as yet unpublished.}. 
However, our most stringent assumption bears on the argument
of $\mes$, as we require the Radon-Nikodym derivative $d\mes/d|\mes|$
to be of bounded variation on $\supp(\mes)$. This assumption, introduced in
\cite{thKus,BKT05}, unlocks many difficulties and will lead us to the
weak convergence of the poles and to the convergence in capacity 
on $\overline\C\setminus(\rspp\cup \supp(\mes))$ of multipoint 
Pad\'e approximants to functions of the form (\ref{fH}).
Moreover we shall prove that each pole of $R$ attracts at least as many
poles of the approximants as its multiplicity, and not much more.
In fact, our hypotheses give rise to an explicit
upper bound on the number of poles of the approximants that may lie
outside a given  neighborhood of
the singular set of $F$. Hence, on each compact subset $K$ of 
$\overline\C\setminus(\rspp\cup \supp(\mes))$, every sequence of approximants
contains a subsequence that converges \emph{uniformly} to $F$ locally
uniformly on $K\setminus E$, where $E$ consists of 
boundedly many (unknown) points. When $\supp(\mes)$ is a finite
union of intervals, results of this type were obtained under stronger
assumptions in \cite{Nut90} for classical Pad\'e approximants.

Finally, we would like to mention that the presented approach can also be carried out for AAK-type meromorphic approximants. Although their definition is rather simple, deriving  functional decomposition for them is not trivial (cf. \cite{AAK71} and \cite{BS02}) and the latter gives rise to more complicated  orthogonality relations than those satisfied by the denominators of Pad\'e approximants. Thus, we consider meromorphic approximants separately in  \cite{uBY1}.

\refstepcounter{section}
\section*{\normalsize\centering\thesection~ Pad\'e Approximation}
\label{sec:pade}

We start by describing the class of measures that we allow in (\ref{fH}) and placing restrictions on the points with respect to which we shall define Pad\'e approximants.

Let $\mes$ be a complex Borel measure whose support
$\spp:=\supp(\mes)\subset\R$ is compact and consists of infinitely many points. 
Denote by $|\mes|$ the total variation measure. 
Clearly $\mes$ is absolutely continuous with respect to $|\mes|$, and we 
shall assume that its Radon-Nikodym derivative (which is of unit
modulus $|\mes|$-a.e.) is of bounded variation. 
In other words, $\mes$ is of the form
\begin{equation}
\label{eq:mesDecomp}
d\mes(t)=e^{i\arm(t)}d|\mes|(t),
\end{equation}
for some real-valued argument function $\arm$ such that\footnote{Note that $e^{i\arm}$ has  bounded variation if and only if $\arm$ can be chosen of bounded  variation.}
\begin{equation}
\label{boundedvarphi}
V(\arm,\spp):=\sup\left\{\sum_{j=1}^N|\arm(x_j)-\arm(x_{j-1})|\right\}<\infty,
\end{equation}
where the supremum is taken over all finite sequences $x_0<x_1<\ldots<x_N$ in $\spp$ as $N$ ranges over $\N$.  

For convenience, we extend the definition of $\arm$ to the whole of $\R$ as follows. Let $I:=[a,b]$ be the convex hull of $\spp$. It is easy to see that if we interpolate $\arm$ linearly in each component of $I\setminus\spp$ and if we set $\arm(x):=\lim_{t\to a, \; t\in\spp}\arm(t)$ for $x<a$ and $\arm(x):=\lim_{t\to b, \; t\in\spp}\arm(t)$ for $x>b$ (the limits exist by (\ref{boundedvarphi})),  the variation of $\arm$ will remain the same. In other words, we may arrange things so that the extension of $\arm$, still denoted by $\arm$, satisfies
\[V(\arm,\spp)=V(\arm,\R)=:V(\arm).\]

Among all complex Borel measures of type  (\ref{eq:mesDecomp})-(\ref{boundedvarphi}), we shall consider only a subclass $\CoM$ defined as follows. We say that a complex measure $\mes$, compactly supported on $\R$, belongs to the class $\CoM$ if it has  an argument of bounded variation and if moreover
\begin{itemize}
\item[(1)] {\it $\supp(\mes)$ is a regular set};
\item[(2)] {\it there exist positive constants $c$ and $L$ such that, for any $x\in\supp(\mes)$ and $\delta\in(0,1)$,  the total variation of $\mes$ satisfies $|\mes|([x-\delta,x+\delta])\geq c\delta^L$}.
\end{itemize} 

In what follows we consider only functions of the form
\begin{equation}
\label{eq:mainFun}
F(z):=\int\frac{d\mes(\xi)}{z-\xi}+R_s(z),
\end{equation}
with $\mes\in\CoM$ and $R_s$ a rational function of type $(s-1,s)$ assumed to be in irreducible form. Hereafter we shall denote by
\begin{equation}
\label{defQs}
Q_s(z)=\prod_{\eta\in\rspp}(z-\eta)^{m(\eta)}
\end{equation}
the denominator of $R_s$, where $\rspp$ is the set of poles of $R_s$ and $m(\eta)$ stands for the multiplicity of $\eta\in\rspp$. Thus, $F$ is a meromorphic function in $\C\setminus\spp$ with poles at each point of $\rspp$ and therefore it is holomorphic in $\overline\C\setminus\wspp$, where
\[\wspp:=\spp\cup\rspp.\]
Note that $F$ does not reduce to a rational function since $\spp$ consists of infinitely many points (cf. \cite[Sec. 5.1]{BMSW06} for a detailed argument).

Diagonal  multipoint  Pad\'e approximants to $F$ are rational functions of type $(n,n)$ that interpolate $F$ at a prescribed system of points. More precisely, pick  $n\in\N$ and let $A_n=\{\zeta_{1,n},\ldots,\zeta_{2n,n}\}$ be a set of $2n$ interpolation points, where the $\zeta_{j,n}\in \overline\C\setminus \wspp$ need not be distinct nor finite. With such an $A_n$ we form the monic polynomial
\begin{equation}
\label{defiweight}
v_{2n}(z)=\prod_{\zeta_{j,n}\in A_n\cap\C}(z-\zeta_{j,n})
\end{equation}
(note that $v_{2n}$ retains only the interpolation points at finite distance thus it needs not have exact degree $2n$).
\smallskip
\newline
{\it Given $F$ of type (\ref{eq:mainFun}) and $A_n$ as above,
the diagonal multipoint Pad\'e approximant to $F$ associated with $A_n$
is the unique rational function $\Pi_n=p_n/q_n$ where 
the polynomials $p_n$ and $q_n$ satisfy:
\begin{itemize}
\item[\textnormal{(i)}] $\deg p_n\leq n$, $\deg q_n\leq n$, and $q_n\not\equiv0$;
\item[\textnormal{(ii)}] $\left(q_n(z)F(z)-p_n(z)\right)/v_{2n}(z)$ is analytic in $\overline\C\setminus\wspp$;
\item[\textnormal{(iii)}] $\left(q_n(z)F(z)-p_n(z)\right)/v_{2n}(z)=O\left(1/z^{n+1}\right)$ as $z\to\infty$.
\end{itemize}}

A multipoint Pad\'e approximant always exists since the conditions for $p_n$ 
and $q_n$ amount to solving a system of $2n+1$ homogeneous linear equations 
with $2n+2$ unknown coefficients, no solution of which can be such that
$q_n\equiv0$ (we may thus assume that $q_n$ is monic); note that (iii)
entails at least one interpolation condition at infinity and therefore $\Pi_n$ is, in fact, of type $(n-1,n)$.

If we let now $\A:=\{A_n\}_{n\in\N}$ be an {\it interpolation scheme}, i.e. a sequence indexed by $n\in\N$ of sets $A_n$ as above, we get a corresponding sequence $\{\Pi_n\}_{n\in\N}$ of diagonal Pad\'e approximants whose asymptotic behavior can be studied when $n$ gets large. Namely, we shall be interested in three types of questions:
\begin{itemize}
\item[(a)] {\it What is the asymptotic distribution of the poles of Pad\'e approximants to $F$?}
\item[(b)] {\it Do some of these poles converge to the polar singularities of $F$?} 
\item[(c)] {\it What can be said about the convergence of such approximants to $F$?}
\end{itemize}

To be able to provide  answers to these questions, we need to place some constraints on interpolation schemes. An interpolation scheme $\A$ is said to be {\it admissible} if
\begin{itemize}
\item[(1)] {\it $\K(\A)$, the set of the limit points of $\A$, is disjoint from $\rspp\cup I$;}
\item[(2)] {\it the counting measures of the points in $A_k$ converge in the weak$^*$ topology to some Borel measure, say $\sigma$, having  finite logarithmic energy;}
\item[(3)] {\it the argument functions of polynomials $v_{2n}$, associated to $\A$ via (\ref{defiweight}), have uniformly bounded derivatives on $I$.}
\end{itemize}

In other words, we call an interpolation scheme admissible if the interpolation points stay away from the poles of $R_s$ and the convex hull of the support of $\mes$, if there exists a Borel measure $\sigma=\sigma(\A)$ supported on $\K(\A)$ such that
\[\sigma_n:=\frac1n \sum_{j=1}^{2n}\delta_{\zeta_{j,n}}\cws\sigma,\]
and if the norms $\|(v_{2n}/|v_{2n}|)^\prime\|_I$ are uniformly bounded with $n$, where $\|\cdot\|_K$ stands for the supremum norm on a set $K$. We call $\sigma$ the {\it asymptotic distribution} of $\A$. Note that $\K(\A)$ is not necessarily compact. If it is not compact,
the finiteness of the logarithmic energy of $\sigma$ is understood as follows.
Since $\K(\A)$ is closed and does not intersect $\wspp$, there exists
$z_0\in\C\setminus\cup_k A_k$ such that $z_0\notin\K(\A)$. Pick such a $z_0$
and set $M_{z_0}(z):=1/(z-z_0)$. Then, all $M_{z_0}(A_k)$ are contained in 
some compact set and their counting measures converge weak$^*$
to $\sigma^\sharp$ such that $\sigma^\sharp(B):=\sigma(M_{z_0}^{-1}(B))$ 
for any Borel set $B\subset\C$. We say that $\A$ is admissible if $\sigma^\sharp$ has finite logarithmic energy. Obviously, this definition does not depend on a particular choice of $z_0$. Further, as a consequence of (3), there exists a finite constant $V_\A$ satisfying
\begin{equation}
\label{eq:VariationPoly}
V(\arg(v_{2n}),I) \leq V_\A \;\;\; \mbox{for any} \;\;\; n\in\N.
\end{equation}
Notice that (3) is satisfied if, for example, all $A_n$ in $\A$ are conjugate-symmetric. More generally, it can be readily verified that (3) amounts to
\begin{equation}
\label{eq:suffast}
\im\left(\int\frac{d\sigma_n(t)}{z-t}\right)=O\left(\frac1n\right)
\end{equation}
uniformly on $I$, which is exactly what we meant in the introduction when saying that the counting measures of interpolation points should converge sufficiently fast.

The four theorems stated below constitute the main results of the paper. For the notions of potential that we use (logarithmic and Green potentials, balayage, equilibrium distributions, capacity and convergence in capacity) the reader may want to consult the appendix.

\begin{thm}
\label{thm:polesPade}
Let $F$ be given by (\ref{eq:mainFun})-(\ref{defQs}) with $\mes\in\CoM$ and let $\{\Pi_n\}_{n\in\N}$ be a sequence of diagonal multipoint Pad\'e approximants to $F$ that corresponds to an admissible interpolation scheme $\A$ with asymptotic distribution $\sigma$. Then the counting measures of the poles of $\Pi_n$ converge in the weak$^*$ sense to $\widehat\sigma$, the balayage of $\sigma$ onto $\spp$.
\end{thm}

We note that the limit distribution of poles of $\Pi_n$ can also be interpreted as the {\it weighted equilibrium distribution} on $\spp$ in the presence of the external field $-U^\sigma$ (cf. \cite[Ch. I]{SaffTotik}).

Recall (cf. \cite[pg. 118]{SaffTotik}) that  $\widehat\delta_\infty$ is simply $\mu_\spp$, the logarithmic equilibrium distribution on $\spp$. Therefore for classical Pad\'e approximants (when each $v_{2n}\equiv1$, i.e. when all the interpolation points are at infinity), the above theorem reduces to the following result.

\begin{cor}
\label{cor:PadeClassique}
Let $F$ be given by (\ref{eq:mainFun})-(\ref{defQs}) with $\mes\in\CoM$ and let $\{\Pi_n\}_{n\in\N}$ be the sequence of Pad\'e approximants to $F$ at infinity. Then the counting measures of the poles of $\Pi_n$ converge to $\mu_\spp$ in the weak$^*$ sense.
\end{cor}

The previous theorem gave one answer to question $(a)$. Our next result addresses question $(c)$ by stating that the approximants behave rather nicely toward the approximated function, namely they converge {\it in 
capacity} to $F$ on $\C\setminus\spp$.

\begin{thm}
\label{thm:convergenceCapacity}
Let $F$, $\A$, and $\{\Pi_n\}_{n\in\N}$ be as in Theorem \ref{thm:polesPade}.  Then 
\begin{equation}
\label{eq:Convergence}
|(F-\Pi_n)(z)|^{1/2n} \cic \exp\left\{-U_{\C\setminus\spp}^\sigma(z)\right\}
\end{equation}
on compact subsets of $\C\setminus\spp$, where  $U_{\C\setminus\spp}^\sigma$ is the Green potential of $\sigma$ relative to $\C\setminus\spp$ and $\cic$ denotes convergence in capacity.
\end{thm}

Finally, we approach question $(b)$. In order to provide an answer to this question, we need some notation. For any point $z\in\C$ define the lower and upper characteristic $\overline m(z),\underline m(z)\in\Z_+$ as
\[\overline m(z) := \inf_U \overline m(z,U), \;\;\; \overline m(z,U) := \lim_{N\to\infty}\max_{n\geq N}\#\{S_n\cap U\},\]
and
\[\underline m(z) := \inf_U \underline m(z,U), \;\;\; \underline m(z,U) := \lim_{N\to\infty}\min_{n\geq N}\#\{S_n\cap U\},\]
respectively, where the infimum is taken over all open sets containing $z$ and $S_n$ is the set of poles of $\Pi_n$, counting multiplicities. Clearly, $\underline m(z)\leq \overline m(z)$, $\underline m(z)=+\infty$ if $z\in\spp$ by Theorem \ref{thm:polesPade}, and $\overline m(z)=0$ if and only if $z$ is not a limit point of poles of $\Pi_n$. Further, let $I_m:=\{[a_j,b_j]\}_{j=1}^m$ be any finite system of intervals covering $\spp$. Also, let $\Arg(\xi)\in(-\pi,\pi]$  be the principal branch of the argument, where we set $\Arg(0)=\pi$. With this definition, $\Arg(\cdot)$ becomes a left continuous function on $\R$. Now, for any interval $[a_j,b_j]$ in $I_m$ we define the angle in which this interval is seen at $\xi\in\C$ by
\[\ang(\xi,[a_j,b_j]):=|\Arg(a_j-\xi)-\Arg(b_j-\xi)|.\]
Finally, we define additively this angle for the whole system, i.e. the angle in which $I_m$ is seen at $\xi$ is defined by\footnote{The notation  does not reflect the dependency on the system of intervals, but the latter will always be made clear.}
\begin{equation}
\label{eq:agnleSys}
\angs(\xi):=\sum_{j=1}^m\ang(\xi,[a_j,b_j]).
\end{equation}
Note that $0\leq\angs(\xi)\leq\pi$ and $\angs(\xi)=\pi$ if and only if $\xi\in I_m$.

The forthcoming theorem implies that each pole of $F$ attracts at least as many poles of Pad\'e approximants as its multiplicity and not much more.

\begin{thm}
\label{thm:tracingPoles}
Let $F$, $\A$, and $\{\Pi_n\}_{n\in\N}$  be as in Theorem \ref{thm:polesPade} and $\angs(\cdot)$ be the angle function for a system of $m$ intervals covering $\spp$. Then
\begin{equation}
\label{eq:LowerM}
\underline m(\eta)\geq m(\eta), \;\;\; \eta\in\rspp,
\end{equation}
and
\begin{equation}
\label{eq:UpperM}
\sum_{\eta\in\rspp\setminus\spp}(\overline m(\eta)-m(\eta))(\pi-\angs(\eta)) \leq V,
\end{equation}
with
\begin{equation}
\label{eq:V}
V := V(\arm) + V_\A+(m+2s^\prime-1)\pi+2\sum_{\eta\in\rspp\setminus\spp}m(\eta)\angs(\eta),
\end{equation}
where $V_\A$ was defined in (\ref{eq:VariationPoly}) and $s^\prime$ is the number of poles of $R$ on $\spp$ counting  multiplicities.  
\end{thm}

The basis of our approach lies in analyzing the asymptotic zero distribution of certain non-Hermitian orthogonal polynomials. It is easy to understand why. Indeed, let $\Gamma$ be any closed Jordan curve that separates $\wspp$ and $\K(\A)$ and contains $\wspp$ in the bounded component of its complement, say $D$. Since
\[(q_nF-p_n)(z)/v_{2n}(z)=O(1/z^{n+1}) \;\;\; \mbox{as} \;\;\; z\to\infty\]
and the left-hand side is analytic in $\C\setminus\wspp$, the Cauchy formula yields
\[\int_\Gamma z^jq_n(z)F(z)\frac{dz}{v_{2n}(z)}=0, \;\;\; j=0,\ldots,n-1, \;\;\; z\in D.\]
Clearly, by writting $R_s$ as
\[
R_s(z) = \sum_{\eta\in\rspp}\sum_{k=0}^{m(\eta)-1}\frac{r_{\eta,k}}{(z-\eta)^{k+1}},
\]
we see that the last equations are equivalent to
\begin{equation}
\label{eq:OrthPade}
\int P_{n-1}(t)q_n(t)\frac{d\mes(t)}{v_{2n}(t)}+\sum_{\eta\in\rspp}\sum_{k=0}^{m(\eta)-1} \frac{r_{\eta,k}}{k!}\left.\left(\frac{P_{n-1}(t)q_n(t)}{v_{2n}(t)}\right)^{(k)}\right|_{t=\eta} = 0
\end{equation} 
for all $P_{n-1}\in\Pc_{n-1}$ by the definition of $F$, the Fubini-Tonelli's theorem, and the residue formula. So, upon taking $P_{n-1}$ to be a multiple of $Q_s$, these relations yield for $n>s$
\begin{equation}
\label{eq:orthRelPade}
\int t^kQ_s(t)q_n(t)\frac{d\mes(t)}{v_{2n}(t)} = 0, \;\;\; k=0,\ldots,n-s-1.
\end{equation}
Hence the denominators of the multipoint Pad\'e approximants to $F$ are polynomials satisfying non-Hermitian orthogonality relations with varying complex measures $d\lambda/v_{2n}$.

The following theorem describes the zero distribution of the polynomials $q_n$ satisfying (\ref{eq:orthRelPade}). Let us stress that, in general, such polynomials need not be unique up to a multiplicative constant nor have exact degree $n$. In the theorem below, it is understood that $q_n$ is {\it any} sequence of such polynomials and that their counting measures are normalized by $1/n$ so that they may no longer be probability measures. This is of no importance since the defect $n-\deg(q_n)$ is uniformly bounded as will be shown later.

\begin{thm}
\label{thm:weakLim}
Let $\{q_n\}_{n\in\N}$ be a sequence of polynomials of degree at most $n$ satisfying weighted orthogonality relations (\ref{eq:orthRelPade}), where $\{v_{2n}\}_{n\in\N}$ is the sequence of monic polynomials associated via (\ref{defiweight}) to some admissible interpolation scheme $\A$ with asymptotic distribution $\sigma$ and where  $\mes\in\CoM$. Then the counting measures $\nu_n$ of the zeros of $q_n(z)=\prod(z-\xi_{j,n})$, namely $\nu_n:=(1/n)\sum\delta_{\xi_{j,n}}$, converge in the weak$^*$ sense to $\widehat\sigma$, the balayage of $\sigma$ onto $\spp=\supp(\mes)$.
\end{thm}

By virtue of the results in the PhD thesis of R. K\"ustner \cite{thKus}, a generalization of the previous theorem can be proved when the  measure $\mes$, instead of belonging to $\CoM$, has an argument of bounded variation and satisfies the 
so-called $\Lambda$-criterion introduced in \cite[Sec. 4.2]{StahlTotik}:
\[\cp\left(\left\{t\in\spp: \; \limsup_{r\to0}\frac{\Log(1/\mu[t-r,t+r])}{\Log(1/r)}<+\infty\right\}\right)=\cp(\spp).\]
However, this assumption would make the exposition heavier and we leave it to the interested reader to carry out the details.

\refstepcounter{section}
\section*{\normalsize\centering\thesection~ Proofs}
\label{sec:proofs}

We start by stating several auxiliary results that are crucial for the proof of Theorem \nolinebreak\ref{thm:weakLim}.
\smallskip
\newline
{\bf Lemma (\cite[Lem. 3.2]{BKT05})}{\it ~Let $\nu$ be a positive measure which has infinitely many points in its support and assume the latter is covered by finitely many disjoint intervals: $\supp(\nu)\subseteq\cup_{j=1}^m[a_j,b_j]$. Let further $\psi$ be a function of bounded variation on $\supp(\nu)$. If the polynomial $u_l(z)=\prod_{j=1}^{d_l}(z-\xi_j)$, $d_l\leq l$, satisfies
\[
\int t^ku_l(t) e^{i\psi(t)}d\nu(t) = 0, \quad k=0,\ldots,l-1,
\]
then
\[
\sum_{j=1}^{d_l}(\pi-\angs(\xi_j)) + (l-d_l)\pi \leq \sum_{j=1}^mV(\psi,[a_j,b_j]) + (m-1)\pi,
\]
where $\angs(\cdot)$ is the angle function defined in (\ref{eq:agnleSys}) for a system of intervals $\cup_{j=1}^m[a_j,b_j]$.}
\smallskip

As a consequence of this lemma, we get the following.
\begin{lem}
\label{lem:aboutAngles}
Let $q_n(z)=\prod_{j=1}^{d_n}(z-\xi_{j,n})$ be an $n$-th orthogonal polynomial in the sense of (\ref{eq:orthRelPade}), where $\mes\in\CoM$ and the polynomials $v_{2n}$ are associated to an admissible interpolation scheme $\A$. Then
\begin{equation}
\label{eq:angleBound}
\sum_{j=1}^{d_n}(\pi-\angs(\xi_{j,n})) + (n-d_n)\pi \leq V(\arm)+V_\A+\sum_{\eta\in\rspp}m(\eta)\angs(\eta)+(m+s-1)\pi,
\end{equation}
where $V_\A$ was defined in (\ref{eq:VariationPoly}) and $\angs(\cdot)$ is the angle function defined in (\ref{eq:agnleSys}) for a system of intervals $I_m:=\cup_{j=1}^m[a_j,b_j]$ that covers $\spp$ with $I=[a_1,b_m]$ being the convex hull of $\spp$.
\end{lem}
{\it Proof:} Denote by $\psi_n(t)$ an argument function for $e^{i\arm(t)}Q_s(t)q_n(t)/v_{2n}(t)$ on $I$, say
\[\psi_n(t)=\arm(t)-\arg(v_{2n}(t))+\sum_{\eta\in\rspp}m(\eta)\Arg(t-\eta)+\sum_{i=1}^{d_n}\Arg(t-\xi_{i,n}).\]
It is easy to see that $\psi_n$ is of bounded variation.  Further, set $l=n-s$,
\[
\psi=\psi_n, \quad \mbox{and}, d\nu(t)=\left|\frac{Q_s(t)q_n(t)}{v_{2n}(t)}\right|d|\mes|(t).
\]
Then it follows from orthogonality relations (\ref{eq:orthRelPade}) that
\[
\int t^k e^{i\psi(t)}d\nu(t) = 0, \quad k=0,\ldots,n-s-1.
\]
Thus, the previous lemma, applied with $u_l\equiv1$, implies that
\[\sum_{j=1}^mV(\psi_n,[a_j,b_j])\geq(n-s-m+1)\pi.\]
So, we are left to show that
\[\sum_{j=1}^mV(\psi_n,[a_j,b_j])\leq V(\arm)+V_\A+\sum_{\eta\in\rspp}m(\eta)\angs(\eta)+\sum_{i=1}^{d_n}\angs(\xi_{i,n}).\]
By the definition of $\psi_n$, we have
\begin{eqnarray}
\sum_{j=1}^mV(\psi_n,[a_j,b_j]) &\leq& \sum_{j=1}^mV(\arm,[a_j,b_j]) + \sum_{j=1}^mV(\arg(v_{2n}),[a_j,b_j]) \nonumber \\ 
{} &+& \sum_{j=1}^m\sum_{\eta\in\rspp}m(\eta)V(\Arg(\cdot-\eta),[a_j,b_j]) \nonumber \\
{} &+& \sum_{j=1}^m\sum_{i=1}^{d_n}V(\Arg(\cdot-\xi_{i,n}),[a_j,b_j]). \nonumber
\end{eqnarray}
The assertion of the lemma now follows from the fact that, by monotonicity, 
\[V(\Arg(\cdot-\xi),[a,b])=\ang(\xi,[a,b]).\]
\boite 

Finally, we state the last two technical observation.
\begin{lem}
\label{lem:aux}
With the previous notation the following statements hold true
\begin{itemize}
\item [(a)] Let $\psi$ be a real function of bounded variation on an interval $[a,b]$ and $Q$ a polynomial. Then there exists a polynomial $T\neq0$ and a constant $\beta\in(0,\pi/32)$ such that
\begin{equation}
\label{eq:polT}
\left|\Arg\left(e^{i\psi(x)}Q(x)T(x)\right)\right|\leq\pi/2-2\beta
\end{equation}
for all $x\in[a,b]$ such that $T(x)Q(x)\neq0$.
\item [(b)] Assume that the polynomials $v_{2n}$ are associated to an admissible interpolation scheme. Then for every $\epsilon>0$ there exists an integer $l$ and a polynomial $T_{l,n}$ of degree at most $l$ satisfying:
\[\left|\frac{v_{2n}(x)}{|v_{2n}(x)|}-T_{l,n}(x)\right|<\epsilon, \;\;\; x\in I,\]
for all $n$ large enough. In particular, the argument of $T_{l,n}/v_{2n}$ lies in the interval $(-2\epsilon,2\epsilon)$ for such $n$.
\end{itemize}
\end{lem}
{\it Proof:} $(a)$ When $Q\equiv1$ this is exactly the statement of Lemma 3.4 in \cite{BKT05} and since  $\psi(x)+\Arg(Q(x))$ is still a real function of bounded variation on $I$, (\ref{eq:polT}) follows.
\smallskip
\newline
$(b)$ This claim follows from Jackson's theorem \cite[Thm. 6.2]{DeVoreLorentz} since the derivatives of $v_{2n}/|v_{2n}|$ are uniformly bounded on $I$.
\boite
\smallskip

Note that Lemma \ref{lem:aboutAngles}, applied with $m=1$, implies that the defect $n-d_n$ is bounded above independently of $n$.

\begin{cor}
\label{cor:numZeros}
Let  $U$ be a neighborhood of $\spp$. Then there exists a constant $k_U\in\N$ such that each $q_n$ has at most $k_U$ zeros outside of $U$ for $n$ large enough.
\end{cor}
{\it Proof:} Since $U$ is open, its intersection with $(-1,1)$ is a countable union of intervals. By compactness, a finite number of them will cover $\spp$, say $\cup_{j=1}^m(a_j,b_j)$. Apply Lemma \ref{lem:aboutAngles} to the closure of these intervals intersected with $I$ and observe that any zero of $q_n$ which lies outside of $U$ will contribute to the left-hand side of (\ref{eq:angleBound}) by more than some positive fixed constant which depends only on $U$. Since the right-hand side of (\ref{eq:angleBound}) does not depend on $n$ and is finite we can have only finitely many such zeros.
\boite
\smallskip
\newline
{\it Proof of Theorem \ref{thm:weakLim}:} Observe that we may  suppose $\A$ is contained in a compact set. 
Indeed, if this is not the case, we can pick a real number 
$x_0\notin\K(\A)\cup\rspp\cup I$ and consider the analytic automorphism of $\overline\C$ given by
$M_{x_0}(z):=1/(z-x_0)$, with inverse 
$M_{x_0}^{-1}(\tau)=x_0+1/\tau$. If we put
$A_n^\sharp:=M_{x_0}(A_n)$, then $\A^\sharp=\{A_n^\sharp\}$ is an admissible 
interpolation scheme having asymptotic distribution $\sigma^\sharp$, 
with $\sigma^\sharp(B)=\sigma(M_{x_0}^{-1}(B))$ 
for any Borel set $B\subset\C$. 
Moreover, the choice of $x_0$ yields that $\K(\A^\sharp)$ is compact. 
Now, if we let
\begin{eqnarray}
\ell_n(\tau)      &=& \tau^nq_n\left(M_{x_0}^{-1}(\tau)\right), \nonumber \\
L_s(\tau)         &=& \tau^sQ_s\left(M_{x_0}^{-1}(\tau)\right), \nonumber \\
P_{n-s-1}^\sharp(\tau) &=& \tau^{n-s-1}P_{n-s-1}\left(M_{x_0}^{-1}(\tau)\right), \nonumber \\
v_{2n}^\sharp(\tau)    &=& \tau^{2n}v_{2n}\left(M_{x_0}^{-1}(\tau)\right), \nonumber
\end{eqnarray}
then $\ell_n$ is a polynomial of degree $n$ with zeros at 
$M_{x_0}(\xi_{j,n})$, $j=1,\ldots,d_n$, and a zero at the origin with 
multiplicity $n-d_n$.  In addition, $v_{2n}^\sharp$ is a polynomial with a 
zero at each point of $A_n^\sharp$, counting multiplicity. Thus, up to a 
multiplicative constant, $v_{2n}^\sharp$ is the polynomial associated 
with $A_n^\sharp$ via (\ref{defiweight}). 
Analogously, $L_s$ is a polynomial of degree $s$ 
with a zero of multiplicity $m(\eta)$  at $M_{x_0}(\eta)$, 
$\eta\in\rspp$, 
and $P^\sharp_{n-s-1}$ is an arbitrary polynomial of degree at most $n-s-1$. 
Making the substitution $t=M_{x_0}^{-1}(\tau)$ in (\ref{eq:orthRelPade}), 
we get
\[\int_{M_{x_0}(\spp)}P^\sharp_{n-s-1}(\tau)L_s(\tau)\ell_n(\tau)\frac{d\mes^\sharp(\tau)}{v_{2n}^\sharp(\tau)}=0, \;\;\; P^\sharp_{n-s-1}\in\Pc_{n-s-1},\]
where $d\mes^\sharp(\tau)=\tau d\mes\left(M_{x_0}^{-1}(\tau)\right)$ is a 
complex measure with compact support $M_{x_0}(\spp)\subset\R$, 
having an argument of bounded variation and total variation measure
$|\mes^\sharp|\in\CoM$. Note that $\tau$ is bounded away from zero on 
$\supp(\mes^\sharp)$, since $\spp$ is compact and therefore bounded away 
from infinity. Now, since Lemma \ref{lem:aboutAngles} implies that 
$n-d_n$ is uniformly bounded above, the asymptotic distribution of 
the counting 
measures of zeros of $\ell_n$ is the same as the asymptotic distribution of 
the images of the counting measures of zeros of $q_n$ under the map $M_{x_0}$. 
As the counting measures of the points in 
$A_n^\sharp$ converge weak$^\sharp$ to $\sigma^\sharp$, it is enough to show that counting measures of zeros of $\ell_n$ converge to $\widehat{\sigma^\sharp}$, since the balayage is preserved 
under $M_{x_0}$ (e.g. because harmonic functions are, {\it cf.} equation (\ref{balayageh}) in the appendix)\footnote{Here we somewhat abuse the notation and use the symbol $\widehat\cdot$ to denote the balayage onto $M_{x_0}(\spp)$, while in the rest of the text it always stands for the balayage onto $\spp$.}. Hence
we assume in the rest of the proof that $\A$ is contained in a compact set, 
say $K_0$, which is disjoint from $\wspp$ by the definition of admissibility.

Now, let $\Gamma$ be a closed Jordan arc such that the bounded component of
$\C\setminus\Gamma$, say $D$, contains $\wspp$ while the unbounded component contains $K_0$. Then $q_n=q_{n,1}\cdot q_{n,2}$, where
\begin{equation}
\label{eq:Laur}
q_{n,1}(z)=\prod_{\xi_{j,n}\in D}(z-\xi_{j,n}) \;\;\; \mbox{and} \;\;\; q_{n,2}(z)=\prod_{\xi_{j,n}\notin D}(z-\xi_{j,n}).
\end{equation}
Corollary \ref{cor:numZeros} assures that degrees of polynomials 
$q_{n,2}$ are uniformly bounded with respect to $n$, therefore the asymptotic 
distribution of the zeros of $q_{n,1}$ coincides with that of $q_n$. 
Denote by $\nu_{n,1}$ the zero counting measure of $q_{n,1}$ 
normalized with $1/n$. Since all $\nu_{n,1}$ are supported on a fixed
compact set, 
Helly's selection theorem and Corollary \ref{cor:numZeros} yield the existence 
of a subsequence $\N_1$ such that $\nu_{n,1}\cws\nu$ for $n\in\N_1$ and some Borel probability
measure $\nu$ supported on $\spp$; remember the defect $n-\deg(q_{n,1})$ is uniformly bounded which is why $\nu$
is a probability measure in spite of the normalization of $q_{n,1}$ with $1/n$. 

Next, we observe it is enough to show that the logarithmic potential of $\nu-\sigma$ is constant q.e. on $\spp$. Indeed, since $\supp(\sigma)$ is disjoint from $\spp$ and subsequently  $U^{-\sigma}$ is harmonic on $\spp$,  $U^\nu$ is bounded q.e. on $\spp$ under this assumption. Hence, by lower semi-continuity of potentials, $U^\nu$ is bounded everywhere on $\spp$ and therefore $\nu$ has finite energy. The latter is sufficient for $\nu$ to be $C$-absolutely continuous\footnote{A Borel measure $\mu$ is called $C$-absolutely continuous if $\mu(E)=0$ for any Borel polar set.}. Moreover, we also get in this case that $U^{\nu-\widehat\sigma}$ is constant q.e. on $\spp$ by (\ref{eq:equalBal}) and, of course, $\widehat\sigma$ is also $C$-absolutely continuous. Thus, $\nu=\widehat\sigma$ by the second unicity theorem \cite[Thm. II.4.6]{SaffTotik}.

Now suppose that  $U^{\nu-\sigma}$ is a constant q.e. on $\spp$. Then there exist nonpolar Borel subsets of $\spp$, say $E_-$ and $E_+$, and two constants $d$ and $\tau>0$ such that
\[U^{\nu-\sigma}(x)\geq d+\tau, \;\;\; x\in E_+, \;\;\; U^{\nu-\sigma}(x)\leq d-2\tau, \;\;\; x\in E_-.\]
Then we claim that there exists $y_0\in\supp(\nu)$ such that
\begin{equation}
\label{eq:point}
U^{\nu-\sigma}(y_0)>d.
\end{equation}
Indeed, otherwise we would have that
\begin{equation}
\label{eq:ineqP}
U^\nu(x)\leq U^\sigma(x)+d, \;\;\; x\in\supp(\nu).
\end{equation}
Then the  principle of domination \cite[Thm. II.3.2]{SaffTotik} would yield that (\ref{eq:ineqP}) is true for all $z\in\C$, but this would contradict the existence of $E_+$.

Since $\K(\A)$ is contained in the complement of $D$, the sequence of potentials $\{U^{\sigma_n}\}_{n\in\N_1}$ converges to $U^\sigma$ locally uniformly in $D$. This implies that for any given sequence of points $\{y_n\}\subset D$ such that $y_n\to y_0$ as $n\to\infty$, $n\in\N_1$, we have
\begin{equation}
\label{eq:sigmaConv}
\lim_{n\to\infty, \; n\in\N_1}U^{\sigma_n}(y_n)=U^\sigma(y_0).
\end{equation}
On the other hand, by applying the principle of descent \cite[Thm. I.6.8]{SaffTotik} for the above sequence $\{y_n\}$, we obtain
\begin{equation}
\label{eq:nuConv}
\liminf_{n\to\infty, \; n\in\N_1} U^{\nu_{n,1}}(y_n)\geq U^\nu(y_0).
\end{equation}
Combining (\ref{eq:point}), (\ref{eq:sigmaConv}), and (\ref{eq:nuConv}) we get
\begin{equation}
\label{eq:nuSigmaConv}
\liminf_{n\to\infty \; n\in\N_1}U^{\nu_{n,1}-\sigma_n}(y_n)\geq U^{\nu-\sigma}(y_0)>d.
\end{equation}
Since $\{y_n\}$ was an arbitrary sequence in $D$ converging to $y_0$, we deduce from (\ref{eq:nuSigmaConv}) that there exists $\rho>0$ such that, for any $y\in[y_0-2\rho,y_0+2\rho]$ and $n\in\N_1$ large enough, the following inequality holds
\begin{equation}
\label{eq:lowerBound}
U^{\nu_{n,1}-\sigma_n}(y)\geq d.
\end{equation}
Clearly 
\begin{equation}
\label{eq:PotToPoly}
U^{\nu_{n,1}-\sigma_n}(y)=\frac{1}{2n}\log\left|\frac{v_{2n}(y)}{q_{n,1}^2(y)}\right|,
\end{equation}
and therefore inequality (\ref{eq:lowerBound}) can be rewritten as
\[\left|\frac{q_{n,1}^2(y)}{v_{2n}(y)}\right|\leq e^{-2nd}, \;\;\; y\in[y_0-2\rho,y_0+2\rho].\]
for all $n\in\N_1$ large enough. We also remark that the same bound holds if $\{q_{n,1}\}$ is replaced by a sequence of monic polynomials, say $\{u_n\}$, of respective degrees $n+o(n)$, whose counting measures normalized by $1/n$ have asymptotic distribution $\nu$. Moreover, in this case
\begin{equation}
\label{eq:aboveBoundPade}
\left|\frac{q_{n,1}(y)u_n(y)}{v_{2n}(y)}\right|\leq e^{-2nd}
\end{equation}
for any $y\in[y_0-2\rho,y_0+2\rho]$ and all $n\in\N_1$ large enough.

In another connection, since $U^{\nu-\sigma}(x)\leq d-2\tau$ on $E_-$, applying the lower envelope theorem \cite[Thm. I.6.9]{SaffTotik} we get
\begin{equation}
\label{eq:lowerEnv}
\liminf_{n\to\infty, \; n\in\N_1}U^{\nu_{n,1}-\sigma_n}(x)=U^{\nu-\sigma}(x)\leq d-2\tau, \;\;\; \mbox{for q.e.} \;\;\; x\in E_-.
\end{equation}
Let $Z$ be a finite system of points from $I$, to be specified later, and denote for simplicity
\[b_n(z) = q_{n,1}^2(z)/v_{2n}(z).\]
Then by \cite{Ank83,Ank_CTP84} there exists $\spp_0\subset\spp$ such that $\spp_0$ is regular$, \cp(E_-\cap\spp_0)>0$ and $\dist(Z,\spp_0)>0$, where $\dist(Z,\spp_0):=\min_{z\in Z}\dist(z,\spp_0)$. Thus, there exists  $x\in E_-\cap\spp_0$ such that 
\[|b_n(x)| \geq  e^{-2n(d-\tau)}, \;\;\; n\in\N_2\subset\N_1,\]
by (\ref{eq:PotToPoly}) and (\ref{eq:lowerEnv}). Let $x_n$ be a point where $|b_n|$ attains its maximum on $\spp_0$, i.e.
\begin{equation}
\label{eq:MaxB}
M_n := \|b_n\|_{\spp_0} = |b_n(x_n)| \geq e^{-2n(d-\tau)}.
\end{equation}
Since $v_{2n}$ has no zeros in $D$, the function $\log|b_n|$ is subharmonic there. Thus, the two-constant theorem \cite[Thm. 4.3.7]{Ransford} on $D\setminus\spp_0$ yields
\[
\log|b_n(z)|\leq \log(M_n)~ \omega_{D\setminus\spp_0}(z,\spp_0) + 2n\log\left(\frac{d(D)}{\dist(\Gamma,K_0)}\right)(1-\omega_{D\setminus\spp_0}(z,\spp_0)),
\]
$z\in D$, where $\omega_{D\setminus\spp_0}$ is the harmonic measure on $D\setminus\spp_0$, $d(D):=\max\{\diam(D),1\},$ and $\diam(D):=\max_{x,y\in D}|x-y|$. Then we get from (\ref{eq:MaxB}) that
\begin{eqnarray}
\label{eq:boundQuotient}
|b_n(z)| &\leq& M_n\left(\frac{1}{M_n}\right)^{1-\omega_{D\setminus\spp_0}(z,\spp_0)}\left(\frac{d(D)}{\dist(\Gamma,K_0)}\right)^{2n(1-\omega_{D\setminus\spp_0}(z,\spp_0))} \nonumber \\
{} &\leq& M_n\exp\left\{2n\Delta(1-\omega_{D\setminus\spp_0}(z,\spp_0))\right\}, \;\;\; z\in D,
\end{eqnarray}
where
\[\Delta := d-\tau+\log(d(D)/\dist(\Gamma,K_0)).\]
Note that $\Delta$ is necessarily positive otherwise $b_n$ would be constant in $D$ by the maximum principle, which is absurd. Moreover, by the regularity of $S_0$, it is known (\cite[Thm. 4.3.4]{Ransford}) that for any $x\in\spp_0$ 
\[\lim_{z\to x}\omega_D(z,\spp_0)=1\]
uniformly with respect to $x\in\spp_0$. Thus, for any $\delta>0$ there exists $r(\delta)<\dist(\spp_0,\Gamma)$ such that for $z$ satisfying $\dist(z,\spp_0)\leq r(\delta)$ we have
\[1-\omega_{D\setminus\spp_0}(z,S_0) \leq \delta/\Delta.\]
This, together with (\ref{eq:boundQuotient}), implies that for fixed $\delta$, to be adjusted later, we have
\[|b_n(z)|\leq M_ne^{2n\delta}, \;\;\; |z-x_n|\leq r(\delta).\]
Note that $b_n$ is analytic in $D$, which, in particular, yields 
\[b^\prime_n(z)=\frac{1}{2\pi i}\int_{|\xi-x_n|=r(\delta)}\frac{b_n(\xi)}{(\xi-z)^2}d\xi, \;\;\; |z-x_n|<r(\delta).\]
Thus, for any $z$ such that $|z-x_n|\leq r(\delta)/2$ we get
\[|b^\prime_n(z)|\leq\frac{1}{2\pi}\cdot\frac{4M_ne^{2n\delta}}{r^2(\delta)}\cdot2\pi r(\delta)=\frac{4M_ne^{2n\delta}}{r(\delta)}.\]
Now, for any $x$ such that 
\begin{equation}
\label{eq:circle}
|x-x_n|\leq\frac{r(\delta)}{8e^{2n\delta}}
\end{equation}
the mean value theorem yields
\[|b_n(x)-b_n(x_n)|\leq\frac{4M_ne^{2n\delta}}{r(\delta)}|x-x_n| \leq \frac{M_n}{2}.\]
Thus, for $x$ satisfying (\ref{eq:circle}) and $n\in\N_2$ we have
\[|b_n(x)|\geq|b_n(x_n)|-|b_n(x)-b_n(x_n)|\geq M_n-\frac{M_n}{2}=\frac{M_n}{2}\]
and by (\ref{eq:MaxB}) and the definition of $b_n$,
\begin{equation}
\label{eq:belowBoundPade}
\left|\frac{q_{n,1}^2(x)}{v_{2n}(x)}\right| \geq \frac12 e^{-2n(d-\tau)}, \;\;\; n\in\N_2.
\end{equation}

Now, Lemma \ref{lem:aux}(a) guarantees that there exist a polynomial $T$ of degree, say $k$, and a number $\beta\in(0,\pi/32)$ such that  
\[\left|\Arg\left(e^{i\arm(t)}Q_s(t)T(t)\right)\right|\leq \frac{\pi}{2}-2\beta,\]
for all $t\in I$ such that $(TQ_s)(t)\neq0$, where $\arm$ is as in (\ref{eq:mesDecomp}). Moreover, for each $n\in\N_2$, we choose $T_{l,n}$ as in Lemma \ref{lem:aux}(b) with $\epsilon=\delta/3$. Since all $T_{l,n}$ are bounded on $I$ by definition and have respective  degrees at most $l$, which does not depend on $n$, there exists $\N_3\subset\N_2$ such that sequence $\{T_{l,n}\}_{n\in\N_3}$ converges uniformly to some polynomial $T_l$ on $I$. In particular, we have that $\deg(T_l)\leq l$ and the argument of $T_l/v_{2n}$ lies in $(-\delta,\delta)$ for $n\in\N_3$ large enough. Denote by $2\alpha$ the smallest even integer strictly greater than $l+k+s$. As soon as $n$ is large enough, since $y_0\in\supp(\nu)$, there exist $\beta_{1,n},\ldots,\beta_{2\alpha,n}$, zeros of $q_{n,1}$, lying in 
\[\left\{z\in\C: \; \dist\left(z,[y_0-\rho,y_0+\rho]\right)\leq \rho\right\},\]
such that
\[\left|\sum_{j=1}^{2\alpha}\Arg\left(\frac{1}{x-\bar\beta_{j,n}}\right)\right|=\left|\sum_{j=1}^{2\alpha}\Arg\left(x-\bar\beta_{j,n}\right)\right| \leq\beta, \;\;\; x\in \R\setminus[y_0-2\rho,y_0+2\rho].\]
Define for $n\in\N_3$ sufficiently large
\[P^*_n(z)=\frac{\overline{q_n(\overline z)}T(z)T_l(z)}{\prod_{j=1}^{2\alpha}(z-\overline \beta_{j,n})}.\]
Then 
\begin{eqnarray}
\left|\Arg\left(\frac{(P^*_nQ_sq_n)(x)e^{i\arm(x)}}{v_{2n}(x)}\right)\right|  &=& \left|\Arg\left(|q_n(x)|^2\prod_{j=1}^{2\alpha}\frac{1}{(x-\bar\beta_{j,n})}~\frac{T_l(x)}{v_{2n}(x)}(TQ_s)(x)e^{i\arm(x)}\right)\right| \nonumber \\
{} &\leq&  \pi/2-\delta, \nonumber
\end{eqnarray}
for all $x\in I\setminus[y_0-2\rho,y_0+2\rho]$ except if $T(x)Q_s(x)=0$, where $\delta$ is chosen so small that $\delta<\beta/2$. This means that for such $x$
\begin{eqnarray}
\label{eq:bb1}
\re\left(\frac{(P^*_nQ_sq_n)(x)e^{i\arm(x)}}{v_{2n}(x)}\right) &\geq& \sin\delta\left|\frac{(P^*_nQ_sq_n)(x)e^{i\arm(x)}}{v_{2n}(x)}\right| \nonumber \\
{} & =& \sin\delta\left|\frac{q_n^2(x)Q_s(x)T(x)T_l(x)}{v_{2n}(x)\prod_{j=1}^{2\alpha}(x-\bar\beta_{j,n})}\right|.
\end{eqnarray}
Now, denote by $m_{n,1}$ the number of zeros of $q_{n,2}$ (defined in (\ref{eq:Laur})) of modulus at least $2\max_{x\in\spp}|x|$, and put $\alpha_n$ for the inverse of their product. Let $m_{n,2}:=\deg(q_{n,2})-m_{n,1}$. Then for $x\in\spp$ we have
\begin{equation}
\label{eq:bb2}
(\dist(\spp,\Gamma))^{m_{n,2}}(1/2)^{m_{n,1}} \leq |\alpha_n q_{n,2}(x)| \leq \left(3\max_{t\in\spp}|t|\right)^{m_{n,2}}(3/2)^{m_{n,1}},
\end{equation}
and since $m_{n,1}+m_{n,2}=\deg(q_{n,2})$ is uniformly bounded with $n$, so is $\{|\alpha_nq_{n,2}|\}$ from above and below on $\spp$.

Finally, if $x\in\spp\setminus[y_0-2\rho,y_0+2\rho]$ satisfies (\ref{eq:circle}), then by (\ref{eq:belowBoundPade})  the  quantity in (\ref{eq:bb1}) is bounded below by
\begin{eqnarray}
{} && |T(x)Q_s(x)|\frac{\sin\delta~\min_{x\in I}|T_l(x)|~\min_{x\in I}|q_{n,2}(x)|^2}{2(\diam(\spp)+2\rho)^{2\alpha}}e^{-2nd+2n\tau} \nonumber \\
{} &=& \frac{c_1}{|\alpha_n^2|}|T(x)Q_s(x)|e^{-2nd+2n\tau}, \nonumber
\end{eqnarray}
where
\[c_1:=\frac{\sin\delta~\min_{x\in I}|T_l(x)|~\min_{x\in I}|\alpha_nq_{n,2}(x)|^2}{2(\diam(\spp)+2\rho)^{2\alpha}}>0\]
by construction of $T_l$ and (\ref{eq:bb2}). Thus,
\begin{eqnarray}
\label{eq:contrBelow}
{} && \re\left(\int_{\spp\setminus[y_0-2\rho,y_0+2\rho]}|\alpha_n|^2P^*_n(t)Q_s(t)q_n(t)\frac{e^{i\arm(t)}}{v_{2n}(t)}d|\mes|(t)\right) \nonumber \\
{} &\geq& \sin\delta\int_{\spp\setminus[y_0-2\rho,y_0+2\rho]}\left|\alpha_n^2P^*_n(t)Q_s(t)q_n(t)\frac{e^{i\arm(t)}}{v_{2n}(t)}\right|d|\mes|(t) \nonumber \\
{} &\geq&  c_1e^{-2nd+2n\tau}\int_{\spp\cap I_n}|T(t)Q_s(t)|d|\mes|(t)\geq c_2e^{-2nd+2n(\tau-L\delta)},
\end{eqnarray}
where $I_n$ is the interval defined by (\ref{eq:circle}). The last inequality is true by the following argument. Recall that $x_n$, the middle point of $I_n$, belongs to $\spp_0$, where $\dist(\spp_0,Z)>0$ and $Z$ is a finite system of points that we choose now to be the zeros of $TQ_s$ on $I$, if any. Then $TQ_s$, which is independent of $n$, is uniformly bounded below on $I_n$ for all $n$ large enough and (\ref{eq:contrBelow}) follows from this, the second requirement in the definition of $\CoM$, and the fact that $I_n$ and $[y_0-2\rho,y_0+2\rho]$ are disjoint for all $n$ large enough. The latter is immediate  if (\ref{eq:aboveBoundPade}) and (\ref{eq:belowBoundPade}) are compared. 

On the other hand, (\ref{eq:aboveBoundPade}), applied with $u_n=P_n^*(z)/\overline{q_{n,2}(\bar z)}$, and (\ref{eq:bb2}) yield that
\begin{equation}
\label{eq:contrAbove}
\left|\int_{[y_0-2\rho,y_0+2\rho]}|\alpha_n|^2P^*_n(t)Q_s(t)q_n(t)\frac{e^{i\arm(t)}}{v_{2n}(t)}d|\mes|(t)\right|\leq c_3e^{-2nd}.
\end{equation}
This completes the proof, since $\delta$ can be taken such  that $\tau-L\delta>0$ and this would contradict orthogonality relations (\ref{eq:orthRelPade}) because, for $n$ large enough, the integral in (\ref{eq:contrBelow}) is much bigger than in (\ref{eq:contrAbove}).
\boite
\smallskip
\newline
{\it Proof of Theorem \ref{thm:polesPade}:} This theorem follows immediately from Theorem \ref{thm:weakLim} and the considerations leading to (\ref{eq:orthRelPade}). 
\boite 
\smallskip

Before we prove Theorem \ref{thm:convergenceCapacity} we shall need one auxiliary lemma.
\begin{lem}
\label{lem:TwoConstant}
Let $D$ be a domain in $\overline\C$ with non-polar boundary, $K^\prime$ be a compact set in $D$, and $\{u_n\}$ be a sequence of subharmonic functions in $D$ such that
\[
u_n(z) \leq M -\epsilon_n, \;\;\; z\in D,
\]
for some constant $M$ and a sequence $\{\epsilon_n\}$ of positive numbers decaying to zero. Further, assume that there exist a compact set $K^\prime$ and positive constants $\epsilon^\prime$ and $\delta^\prime$, independent of $n$, for which holds
\[
u_n(z) \leq M - \epsilon^\prime, \;\;\; z\in K_n\subset K^\prime, \;\;\; \cp(K_n)\geq\delta^\prime.
\]
Then for any compact set $K\subset D\setminus K^\prime$ there exists a positive constant $\epsilon(K)$ such that
\[
u_n(z) \leq M - \epsilon(K), \;\;\; z\in K,
\]
for all $n$ large enough. 
\end{lem}
{\it Proof:} Let $\omega_n$ be the harmonic measure for $D_n := D\setminus K_n$. Then the two-constant theorem \cite[Thm. 4.3.7]{Ransford} yields that
\begin{eqnarray}
u_n(z) &\leq& (M - \epsilon^\prime)\omega_n(z,K_n) + (M-\epsilon_n)(1-\omega_n(z,K_n)) \nonumber \\
{}     &\leq& M - (\epsilon^\prime-\epsilon_n)\omega_n(z,K_n), \;\;\; z\in D_n. \nonumber
\end{eqnarray}
Thus, we need to show that for any $K\subset D\setminus K^\prime$ there exists a constant $\delta(K)>0$ such that
\[
\omega_n(z,K_n) \geq \delta(K), \;\;\; z\in K.
\]
Assume to the contrary that there exists a sequence of points $\{z_n\}_{n\in\N_1}\subset K$, $\N_1\subset\N$, such that
\begin{equation}
\label{eq:contrary}
\omega_n(z_n,K_n) \to 0 \;\; \mbox{as} \;\; n\to\infty, \;\; n\in\N_1.
\end{equation}
By \cite[Theorem 4.3.4]{Ransford}, $\omega_n(\cdot,K_n)$ is the unique bounded harmonic function in $D_n$ such that
\[
\lim_{z\to\zeta} \omega(z,K_n) = 1_{K_n}(\zeta)
\]
for any regular $\zeta\in\partial D_n$, where $1_{K_n}$ is the characteristic function of $\partial K_n$. Then it follows from (\ref{eq:GreenEqual}) of the appendix that
\begin{equation}
\label{eq:HarmEqPot}
\cp(K_n,\partial D) U_D^{\mu_{(K_n,\partial D)}} \equiv \omega_n(\cdot,K_n),
\end{equation}
where $\mu_{(K_n,\partial D)}$ is the Green equilibrium measure on $K$ relative to $D$. Since all the measures $\mu_{(K_n,\partial D)}$ are supported in the compact set $K^\prime$, there exists a probability measure $\mu$ such that
\[
\mu_{(K_n,\partial D)} \cws \mu \;\; \mbox{as} \;\; n\to\infty, \;\; n\in\N_2\subset\N_1.
\]
Without loss of generality we may suppose that $z_n\to z^*\in K$ as $n\to\infty$, $n\in\N_2$. Let, as usual, $g_D(\cdot,t)$ be the Green function for $D$ with pole at $t\in D$. Then, by the uniform equicontinuity of $\{g_D(\cdot,t)\}_{t\in K^\prime}$ on $K$, we get
\[
U_D^{\mu_{(K_n,\partial D)}}(z_n) \to U_D^{\mu}(z^*) \neq 0 \;\; \mbox{as} \;\; n\to\infty, \;\; n\in\N_2.
\]
Therefore, (\ref{eq:contrary}) and (\ref{eq:HarmEqPot}) necessarily mean that
\begin{equation}
\label{eq:contrary2}
\cp(K_n,\partial D) \to 0 \;\; \mbox{as} \;\; n\to\infty, \;\; n\in\N_2.
\end{equation}
By definition, $1/\cp(K_n,\partial D)$ is the minimum among Green energies of probability measures supported on $K_n$. Thus, the sequence of Green energies of the logarithmic equilibrium measures on $K_n$, $\mu_{K_n}$, diverges to infinity by (\ref{eq:contrary2}). Moreover, since 
\[
\left\{g(\cdot,t)+\log|\cdot-t|\right\}_{t\in K^\prime}
\]
is a family of harmonic functions in $D$ whose moduli are uniformly bounded above on $K^\prime$, the logarithmic energies of $\mu_{K_n}$ diverge to infinity. In other words,
\[
\cp(K_n) \to 0 \;\; \mbox{as} \;\; n\to\infty, \;\; n\in\N_2,
\]
which is impossible by the initial assumptions. This proves the lemma.
\boite
\smallskip
\newline
{\it Proof of Theorem \ref{thm:convergenceCapacity}:} Exactly as in the proof of Theorem \ref{thm:weakLim} we can suppose that all the interpolation points are contained in some compact set $K_0$ disjoint from $\wspp$. By virtue of the  Hermite interpolation formula ( cf. \cite[Lemma 6.1.2, (1.23)]{StahlTotik}), the error $e_n:=F-\Pi_n$ has the following representation
\begin{equation}
\label{eq:form1Pade}
e_n(z)=\frac{v_{2n}(z)}{(p_{n-s}Q_sq_n)(z)}\int \frac{(p_{n-s}Q_sq_n)(t)}{v_{2n}(t)}\frac{d\mes(t)}{z-t}, \;\;\; z\in\C\setminus\spp,
\end{equation}
where $p_{n-s}$ is an arbitrary polynomials in $\Pc_{n-s}$. Since almost all of the zeros of $q_n$ approach $\spp$ by Corollary \ref{cor:numZeros}, we always can fix $s$ of them, say $\xi_{1,n},\ldots,\xi_{s,n}$, in such a manner that the absolute value of
\[
l_{s,n}(z):=\prod_{j=1}^s(z-\xi_{j,n})
\]
is uniformly bounded above and below on any given compact subset $K\subset\C\setminus\spp$ for all $n$ large enough (depending on $K$). In what follows we choose $p_{n-s}(z) := \overline{q_n(\bar z)}/\overline{l_{s,n}(\bar z)}$. Set also $A_n$ to be
\[
A_n(z) := \int \frac{(|p_{n-s}^2|Q_sl_{s,n})(t)}{v_{2n}(t)}\frac{d\mes(t)}{z-t}, \;\;\; z\in\C\setminus\spp.
\]
First, we show that 
\begin{equation}
\label{eq:ConvIn}
|A_n|^{1/2n}\cic\exp\{-c(\sigma,\C\setminus\spp)\}
\end{equation}
on compact subsets of $\C\setminus\spp$, where $c(\sigma,\C\setminus\spp)$ is defined in (\ref{eq:equalBal}) of the appendix. Clearly, for any compact set $K\subset\C\setminus\spp$ there exists a constant $c(K)$, independent of $n$, such that
\begin{equation}
\label{eq:UpperBoundIn}
|A_n(z)| \leq c(K)\left\|\frac{p_{n-s}^2}{v_{2n}}\right\|_\spp, \;\;\; z\in K.
\end{equation}
Let $\nu_n$ and $\sigma_n$ be the counting measures of zeros of $p_{n-s}^2$ and $v_{2n}$, respectively. Then
\begin{eqnarray}
\limsup_{n\to\infty}\left|\frac{p_{n-s}^2(t)}{v_{2n}(t)}\right|^{1/2n} &=& \liminf \exp\left\{U^{\sigma_n-\nu_n}(t)\right\} = \exp\left\{U^{\sigma-\widehat\sigma}(t)\right\} \nonumber \\
\label{eq:LET}
{} &=& \exp\{-c(\sigma;\C\setminus\spp)\} \;\;\; \mbox{q.e. on} \;\;\; \spp
\end{eqnarray}
by Theorem \ref{thm:polesPade}, the lower envelope theorem \cite[Thm. I.6.9]{SaffTotik}, and (\ref{eq:equalBal}) of the appendix. Moreover, by the principle of descent \cite[Thm. I.6.8]{SaffTotik}, we get that
\begin{equation}
\label{eq:POD}
\limsup_{n\to\infty}\left|\frac{p_{n-s}^2(t)}{v_{2n}(t)}\right|^{1/2n} \leq \exp\left\{U^{\sigma-\widehat\sigma}(t)\right\} = \exp\{-c(\sigma;\C\setminus\spp)\}
\end{equation}
uniformly on $\spp$, where the last equality holds by the regularity of $\spp$. Now it is immediate from (\ref{eq:LET}) and (\ref{eq:POD}) that
\begin{equation}
\label{eq:LimitNorms}
\lim_{n\to\infty}\left\|\frac{p_{n-s}^2}{v_{2n}}\right\|_\spp^{1/2n} = \exp\{-c(\sigma;\C\setminus\spp)\}.
\end{equation}
Indeed, since the whole sequence $\{\nu_n\}$ converges to $\widehat\sigma$, (\ref{eq:LET}) holds for any subsequence of $\N$. Thus, there cannot exist a subsequence of natural numbers for which the limit in (\ref{eq:LimitNorms}) would not hold. Suppose now that (\ref{eq:ConvIn}) is false. Then there would exist a compact set $K^\prime\subset\C\setminus\spp$ and $\epsilon^\prime>0$ such that
\begin{equation}
\label{eq:ContrIn}
\cp\left\{z\in K^\prime:~\left||A_n(z)|^{1/2n}-\exp\{-c(\sigma;\C\setminus\spp)\}\right|\geq\epsilon^\prime\right\} \not\to 0.
\end{equation}
Combining (\ref{eq:ContrIn}), (\ref{eq:LimitNorms}), and (\ref{eq:UpperBoundIn}) we see that there would exist a sequence of compact sets $K_n \subset K^\prime$,$\cp(K_n)\geq\delta^\prime>0$, such that
\begin{equation}
\label{eq:FalseUpperBoundIn}
|A_n(z)|^{1/2n} \leq \exp\{-c(\sigma;\C\setminus\spp)\} - \epsilon^\prime, \;\;\; z\in K_n.
\end{equation}
Now, let $\Gamma$ be a closed Jordan curve that separates $\spp$ from $K_0$, a compact set containing all the zeros of $v_{2n}$, and $K^\prime$. Assume further that $\spp$ belongs to the bounded component of the complement of $\Gamma$. Observe that $(1/2n)\log|A_n|$ is a subharmonic function in $\C\setminus\spp$. Then by (\ref{eq:UpperBoundIn}), (\ref{eq:LimitNorms}), and (\ref{eq:FalseUpperBoundIn}) enable us to apply Lemma \ref{lem:TwoConstant} with $M=-c(\sigma;\C\setminus\spp)$ which yields that there exists $\epsilon(\Gamma)>0$ such that
\begin{equation}
\label{eq:FalseOnGamma}
|A_n(z)|^{1/2n} \leq \exp\{-c(\sigma;\C\setminus\spp)-\epsilon(\Gamma)\}
\end{equation}
uniformly on $\Gamma$ and for all $n$ large enough. Define
\[
J_n : = \left|\int_\Gamma T_l(z)T(z)\overline{l_{s,n}(\bar z)}A_n(z)\frac{dz}{2\pi i}\right|,
\]
where the polynomials $T_l$ and $T$ are chosen as in Theorem \ref{thm:polesPade} (see discussion after (\ref{eq:belowBoundPade})). We get from (\ref{eq:FalseOnGamma}) that
\begin{equation}
\label{eq:ContrPart1}
\limsup_{n\to\infty} J_n^{1/2n} \leq \exp\{-c(\sigma;\C\setminus\spp)-\epsilon(\Gamma)\}.
\end{equation}
In another connection, Fubini-Tonelli theorem and the Cauchy integral formula yield
\begin{eqnarray}
J_n &=& \left|\int_\Gamma T_l(z)T(z)\overline{l_{s,n}(\bar z)}\left(\int\frac{(|p_{n-s}^2|Q_sl_{s,n})(t)}{v_{2n}(t)}\frac{d\mes(t)}{z-t}\right) \frac{dz}{2\pi i}\right| \nonumber \\
\label{eq:BelowOnGamma1}
{} &=& \left|\int|q_n^2(t)|\frac{T_l(t)}{v_{2n}(t)}(TQ_s)(t)e^{i\arm(t)}d|\mes|(t)\right|. 
\end{eqnarray}
Exactly as in (\ref{eq:bb1}), we can write
\begin{equation}
\label{eq:BelowOnGamma2}
\re\left(\frac{(T_lTQ_s)(t)e^{i\arm(t)}}{v_{2n}(t)}\right) \geq \sin(\delta)\left|\frac{(T_lTQ_s)(t)}{v_{2n}(t)}\right|, \;\;\; t\in I,
\end{equation}
where $I$ is the convex hull of $\spp$ and $\delta>0$ has the same meaning as in the proof of Theorem \ref{thm:polesPade} (see construction after (\ref{eq:boundQuotient})). Thus, we derive from (\ref{eq:BelowOnGamma1}) and (\ref{eq:BelowOnGamma2}) that
\begin{equation}
\label{eq:BelowOnGamma3}
J_n \geq \sin(\delta) \int |b_n(t)| ~ |(T_lTQ_s)(t)| d|\mes|(t), \;\;\; b_n := q_n^2/v_{2n}.
\end{equation}
Let $\spp_0$ be a closed subset of $\spp$ of positive capacity that lies at positive distance from the zeros of $TQ_s$ on $I$ (as in the proof of Theorem \ref{thm:polesPade}, we refer to \cite{Ank83,Ank_CTP84} for the existence of this set). Further, let  $x_n\in\spp_0$ be such that
\[
\|b_n\|_{\spp_0} = |b_n(x_n)|.
\]
Then it follows from (\ref{eq:LET}) that
\[
\|b_n\|_{\spp_0} \geq \exp\{-2n(c(\sigma;\C\setminus\spp)+\epsilon)\}
\]
for any $\epsilon>0$ and all $n$ large enough. Proceeding as in Theorem \ref{thm:polesPade} (see equations (\ref{eq:circle}) and (\ref{eq:belowBoundPade})), we get that
\begin{equation}
\label{eq:BelowOnGamma4}
|b_n(t)| \geq \frac12\exp\{-2n(c(\sigma;\C\setminus\spp)+\epsilon)\}, \;\;\; t\in I_n,
\end{equation}
where
\[
I_n:=\left\{x\in\spp_0:~|x-x_n|\leq r_\delta e^{-2n\delta}\right\}
\]
and $r_\delta$ is some function of $\delta$ continuous and vanishing at zero. Then by combining (\ref{eq:BelowOnGamma3}) and (\ref{eq:BelowOnGamma4}), we obtain exactly as in (\ref{eq:contrBelow}) that there exists a constant $c_1$ independent of $n$ such that
\[
J_n \geq \sin(\delta) \int_{I_n} |b_n(t)| ~ |(T_lTQ_s)(t)| d|\mes|(t) \geq c_1 \exp\{-2n(c(\sigma;\C\setminus\spp)+\epsilon+L\delta)\}.
\]
Thus, we have that
\begin{equation}
\label{eq:ContrPart2}
\liminf_{n\to\infty} J_n^{1/2n} \geq \exp\{-c(\sigma;\C\setminus\spp)-\epsilon-L\delta\}.
\end{equation}
Now, by choosing $\epsilon$ and $\delta$ small enough that $\epsilon+L\delta < \epsilon(\Gamma)$, we arrive at contradiction between (\ref{eq:ContrPart1}) and (\ref{eq:ContrPart2}). Therefore, the convergence in (\ref{eq:ConvIn}) holds.

Second, we show that
\begin{equation}
\label{eq:ConvPolies}
\left|\frac{v_{2n}(z)l_{s,n}(z)}{q_n^2(z)Q_s(z)}\right|^{1/2n}\cic\exp\left\{c(\sigma;\C\setminus\spp)-U_{\C\setminus\spp}^\sigma(z)\right\}
\end{equation}
on compact subsets of $\C\setminus\spp$. Let $K\subset\C\setminus\spp$ be compact and let $U$ be a bounded open set containing $K$ and not intersecting $\spp$. Define
\[
q_{n,1}(z) := \prod_{\xi\in U:~q_n(\xi)=0}(z-\xi) \;\;\; \mbox{and} \;\;\; q_{n,2}(z) := q_n(z)/q_{n,1}(z).
\]
Corollary \ref{cor:numZeros} yields that there exists fixed $m\in\N$ such that $\deg(q_{n,1})\leq m$. Then
\[
\left|\frac{v_{2n}(z)}{q_{n,2}^2(z)}\right|^{1/2n} \to \exp\left\{U^{\widehat\sigma-\sigma}(z)\right\} = \exp\left\{c(\sigma;\C\setminus\spp)-U_{\C\setminus\spp}^\sigma(z)\right\}
\]
uniformly on $K$ by Theorem \ref{thm:polesPade}, definition of $v_{2n}$, and (\ref{eq:toRemind}) of the appendix. Moreover, it is an immediate consequence the choice of $l_{s,n}$, the uniform boundedness of the degrees of $q_{n,1}^2Q_s$, and \cite[Thm. 5.2.5]{Ransford} ($\cp(\{z:|(q_{n,1}^2Q_s)(z)|\leq\epsilon\})=\epsilon^{\deg(q_{n,1}^2Q_s)}$) that
\[
\left|\frac{l_{s,n}(z)}{q_{n,1}^2(z)Q_s(z)}\right|^{1/2n}\cic1, \;\;\; z\in K.
\]
Thus, we obtain (\ref{eq:ConvPolies}). It is clear now that (\ref{eq:Convergence}) follows from (\ref{eq:form1Pade}), (\ref{eq:ConvIn}), and (\ref{eq:ConvPolies}).
\boite
\smallskip
\newline
{\it Proof of Theorem \ref{thm:tracingPoles}:} Inequality (\ref{eq:LowerM}) is trivial for any $\eta\in\rspp\cap\spp$. Suppose now that $\eta\in\rspp\setminus\spp$ and that $\underline m(\eta)< m(\eta)$. This would mean that there exists an open set $U$,  $U\cap\wspp=\{\eta\}$, such that $\underline m(\eta,U)<m(\eta)$ and therefore would exist a subsequence $\N_1\subset\N$ such that
\[\#\{\spp_n\cap U\}<m(\eta), \;\;\; n\in\N_1.\]
It was proved in Theorem \ref{thm:convergenceCapacity} that $\{\Pi_n\}$ converges in capacity on compact subsets of $\C\setminus\spp$ to $F$. Thus, $\{\Pi_n\}_{n\in\N_1}$ is a sequence of meromorphic (in fact, rational) functions in $U$ with at most $m(\eta)$ poles there, which converges in capacity on $U$ to a meromorphic function $\left.F\right|_{U}$ with exactly one pole of multiplicity $m(\eta)$. Then by Gonchar's lemma \cite[Lemma 1]{Gon75b} each $\Pi_n$ has exactly $m(\eta)$ poles in $U$ and these poles converge to $\eta$. This finishes the proof of (\ref{eq:LowerM}).

Now, for any $\eta\in\rspp\setminus\spp$ the upper characteristic $\overline m(\eta)$ is finite by Corollary \ref{cor:numZeros}. Therefore there exist domains $D_\eta$, $D_\eta\cap\wspp=\{\eta\}$, such that $\overline m(\eta)=\overline m(\eta,D_\eta)$, $\eta\in\rspp\setminus\spp$. Further, let $\angs(\cdot)$ be the angle function defined in (\ref{eq:agnleSys}) for a system of $m$ intervals covering $\spp$ and let $S_n=\{\xi_{1,n},\ldots,\xi_{d_n,n}\}$. Then by Lemma \ref{lem:aboutAngles} we have
\begin{equation}
\label{eq:nameless}
\sum_{j=1}^{d_n}(\pi-\angs(\xi_{j,n})) + (n-d_n)\pi \leq V(\arm)+V_\A+(m+s-1)\pi+\sum_{\zeta\in\rspp}m(\zeta)\angs(\zeta).
\end{equation}
Then for $n$ large enough (\ref{eq:nameless}) yields
\[\sum_{\eta\in\rspp\setminus\spp}\left(\sum_{\xi_{j,n}\in D_\eta}(\pi -\angs(\xi_{j,n}))-m(\eta)(\pi-\angs(\eta))\right) \leq V,\]
where $V$ was defined in (\ref{eq:V}). Thus,
\begin{equation}
\label{eq:nmless}
\begin{array}{rl}
\sum_{\eta\in\rspp\setminus\spp} & \left(\#\{S_n\cap D_\eta\} -m(\eta)\right)(\pi-\angs(\eta)) \smallskip \\
& \leq \sum_{\eta\in\rspp\setminus\spp}\#\{S_n\cap D_\eta\}\left(\max_{\xi\in D_\eta}\angs(\xi)-\angs(\eta)\right) + V
\end{array}
\end{equation}
for all $n$ large enough. However, since $\{\max_{n\geq N}\#\{S_n\cap D_\eta\}\}_{N\in\N}$ is a decreasing sequence of integers, $\overline m(\eta) = \overline m(\eta,D_\eta)= \#\{S_n\cap D_\eta\}$ for infinitely many $n\in\N$. Therefore, we get from (\ref{eq:nmless}) that
\begin{equation}
\label{eq:mnless}
\sum_{\eta\in\rspp\setminus\spp}\left(\overline m(\eta) -m(\eta)\right)(\pi-\angs(\eta)) \leq V + \sum_{\eta\in\rspp\setminus\spp}\overline m(\eta)\left(\max_{\xi\in D_\eta}\angs(\xi)-\angs(\eta)\right).
\end{equation}
Observe now that the left-hand side and the first summand on the right-hand side of (\ref{eq:mnless}) are simply constants. Moreover, the second summand on the right-hand side of (\ref{eq:mnless}) can be maid arbitrarily small by taking smaller neighborhoods $D_\eta$. Thus, (\ref{eq:UpperM}) follows.
\boite

\refstepcounter{section}
\section*{\normalsize\centering\thesection~ Numerical Experiments}
\label{sec:compute}

We restricted ourselves to the case of classical Pad\'e approximants and we constructed their denominators by solving the orthogonality relations (\ref{eq:orthRelPade}) with $w_{2n}\equiv1$. Thus, finding these denominators amounts to solving a system of linear equations whose coefficients are obtained from the moments of the measure $\lambda$.

In the numerical experiments below we approximate 
function $F$ given by the formula
\begin{eqnarray}
F(z) &=& 7\int_{[-6/7,-1/8]}\frac{e^{it}dt}{z-t} - (3+i)\int_{[2/5,1/2]}\frac{t-3/5}{t-2i}\frac{dt}{z-t}  \nonumber \\
&& + (2-4i)\int_{[2/3,7/8]}\frac{\ln(t)dt}{z-t} \nonumber \\ 
&& + \frac{1}{(z+3/7-4i/7)^2} + \frac{2}{(z-5/9-3i/4)^3} + \frac{6}{(z+1/5+6i/7)^4}. \nonumber
\end{eqnarray}

\begin{figure}[h!]
\centering
\includegraphics[scale=.4]{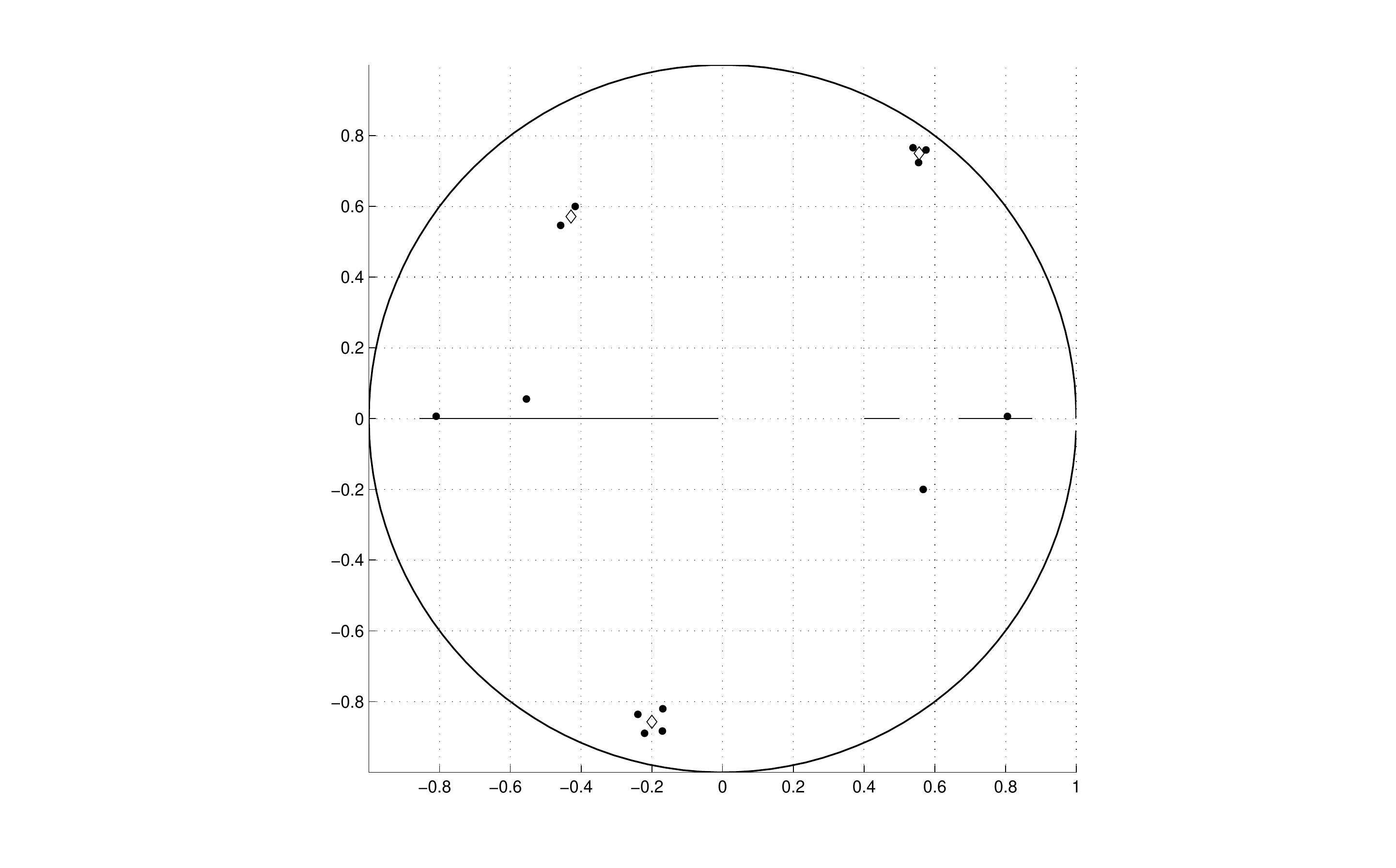}
\includegraphics[scale=.4]{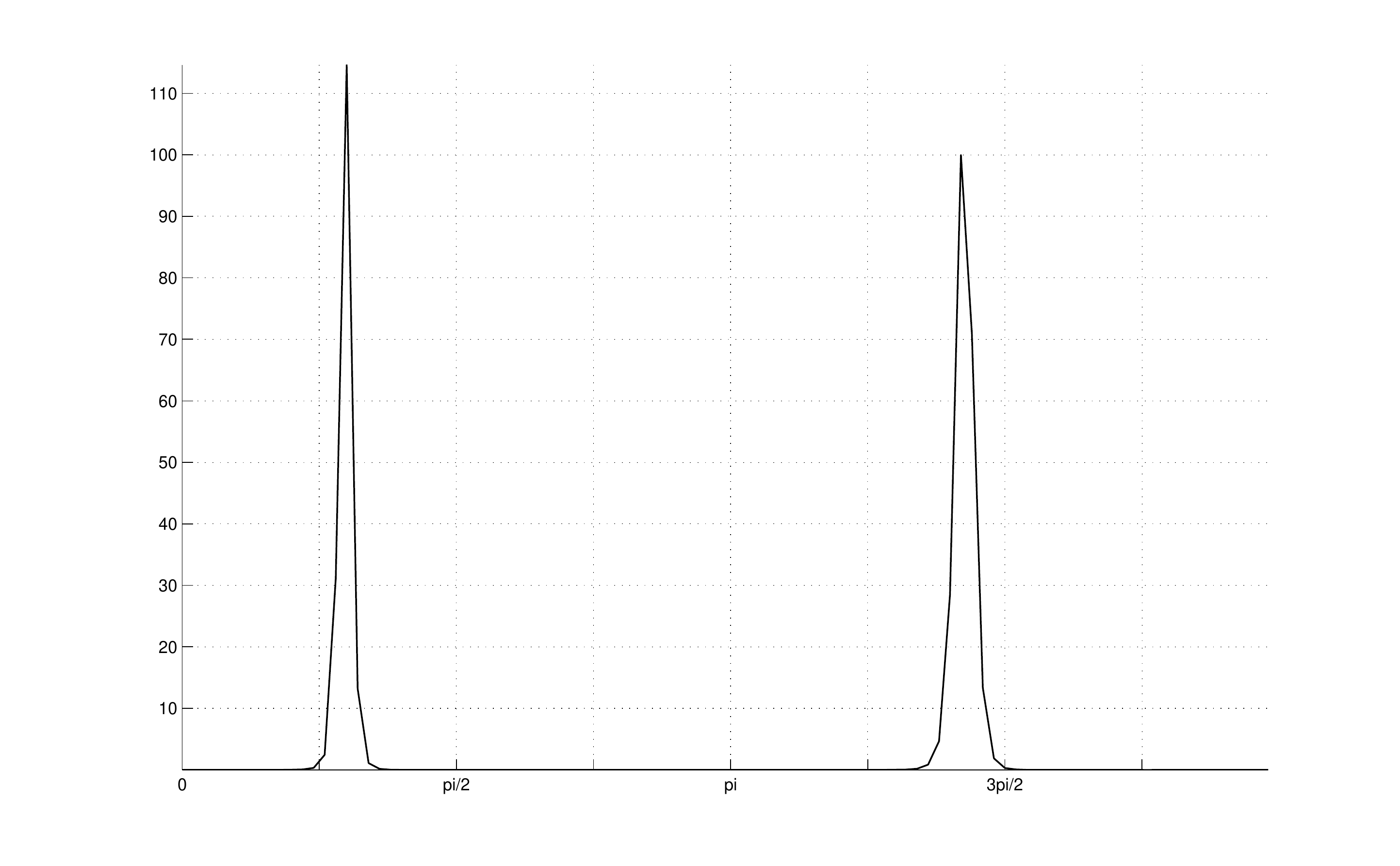}
\caption{\small Poles of $\Pi_{13}$ and the error $|F-\Pi_{13}|$ on $\T$}
\end{figure}
\begin{figure}[h!]
\centering
\includegraphics[scale=.4]{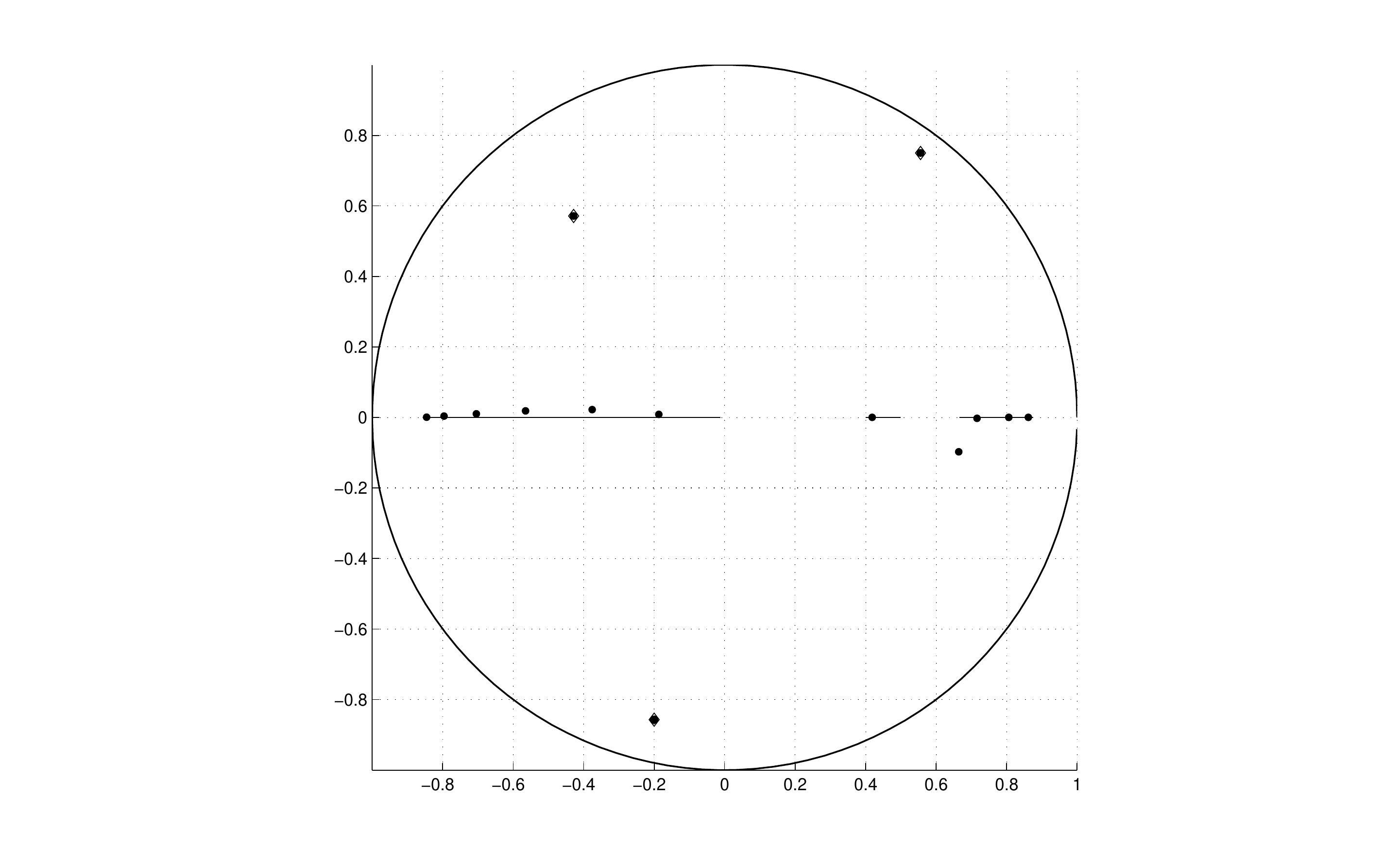}
\includegraphics[scale=.4]{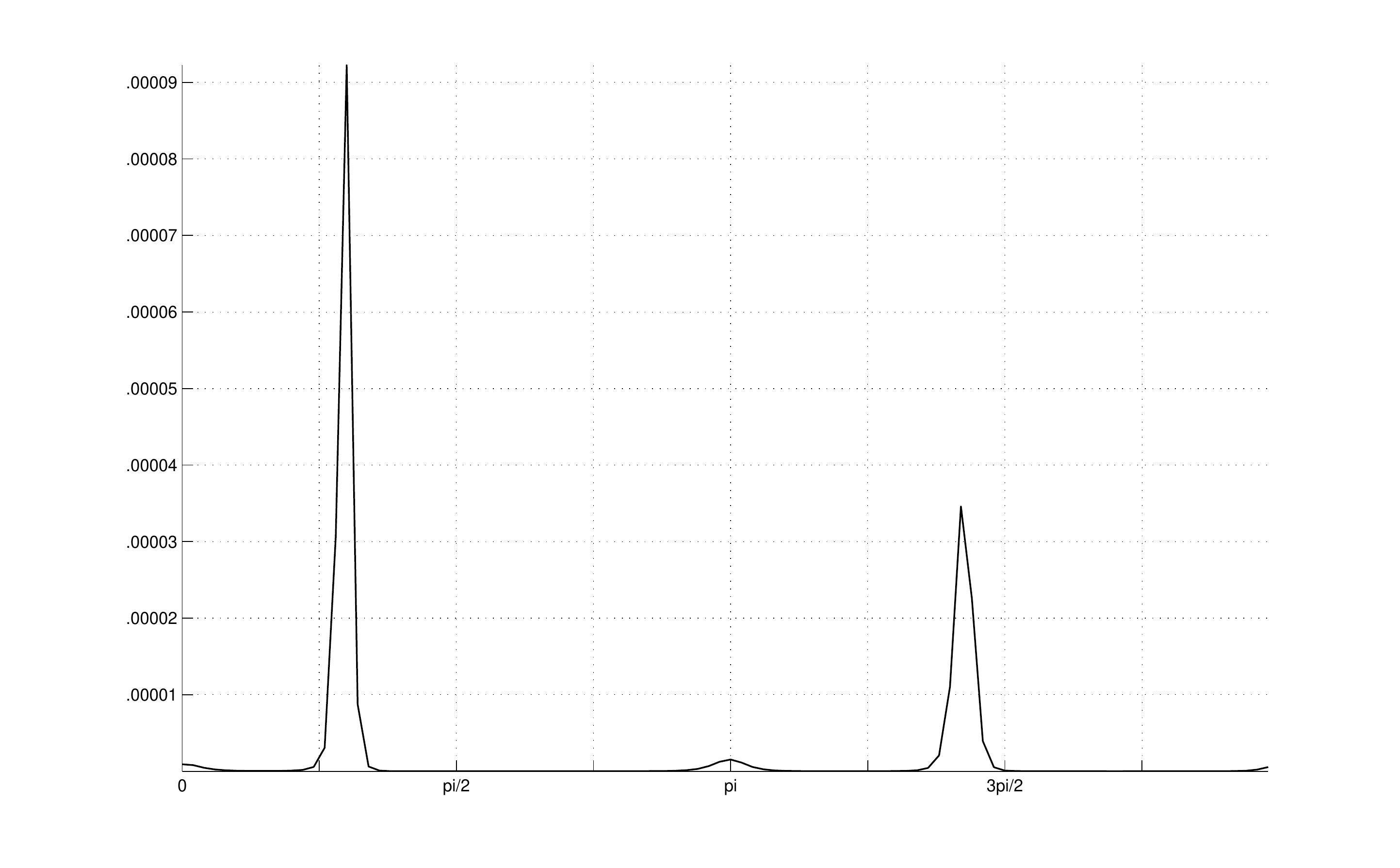}
\caption{\small Poles of $\Pi_{20}$ the error $|F-\Pi_{20}|$ on $\T$}
\end{figure}

On Figures 1a and 2a the solid lines stand for the support of the measure, diamonds depict the polar singularities of $F$, and disks denote the poles of the corresponding approximants. Note that the poles of $F$ seem to attract the singularities first. On Figures 1b and 2b the absolute value of the error on the unit circle is displayed for the corresponding approximants. The horizontal parts of the curves are of magnitude about $10^{-3}$ on Figure 1b and of magnitude about $10^{-9}$ on Figure 2b.

\refstepcounter{section}
\section*{\normalsize\centering\thesection~ Appendix}
\label{prelnot}
\setcounter{section}{1}
\renewcommand{\theequation}{\Alph{section}.\arabic{equation}}

Below we sketch some basic notions of logarithmic potential theory that were used throughout the paper. We refer the reader  to the monographs \cite{Ransford,SaffTotik} for a complete treatment.

The {\it logarithmic potential} and the {\it logarithmic energy} of a finite positive 
measure $\mu$, compactly supported in $\C$, are defined by
\begin{equation}
\label{eq1sL}
U^\mu(z):=\int\log\frac{1}{|z-t|}d\mu(t),~~~~z\in\C,
\end{equation}
and
\begin{equation}
\label{eq2sL}
I[\mu]:=\int U^\mu(z)d\mu(z)=\int\int\log\frac{1}{|z-t|}d\mu(t)d\mu(z),
\end{equation}
respectively. The function $U^\mu$ is superharmonic with values in $(-\infty,+\infty]$, and is not identically $+\infty$. It is bounded below on $\supp(\mu)$ so that
$I[\mu]\in(-\infty,+\infty]$.

Let now $E\subset \C$ be compact and $\Lm(E)$ denote the set of all 
probability measures supported on $E$. If the logarithmic energy of every 
measure in $\Lm(E)$ is infinite, we say that $E$ is {\it polar}. 
Otherwise, there exists a unique $\mu_E\in\Lm(E)$ that
minimizes the logarithmic energy over all measures in $\Lm(E)$. 
This measure is called the {\it equilibrium distribution} on $E$.  The {\it logarithmic capacity}, or simply
the capacity, of $E$ is defined as
\[\cp(E)=\exp\{-I[\mu_E]\}.\]
By definition, the capacity of an arbitrary subset of $\C$ is the {\it supremum} of the capacities of its compact subsets. We agree that the capacity of a polar set is zero. We define convergence in capacity as follows. We say that a sequence of functions $\{h_n\}$ {\it converges in capacity} to a function $h$ on a compact set $K$ if for any $\epsilon>0$ holds
\[
\cp\left(\{z\in K:~|(h_n-h)(z)|\geq\epsilon\}\right)\to0 \;\;\; \mbox{as} \;\;\; n\to\infty.
\]
We also say that a sequence converges in capacity in an open set $\Omega$ if it is converges in capacity on any compact subset of $\Omega$.

Another important concept is the \emph{regularity} of a compact set. We restrict to the case when $E$ has connected complement, say $\Omega$. Then $E$ is called regular if the Dirichlet problem on $\partial \Omega$ is solvable, in other words, if any continuous function on $\partial \Omega$ is the trace (limiting boundary values) of some function harmonic in $\Omega$. Thus, regularity is a property of $\partial\Omega$ rather than $E$ itself. It is also  known \cite[pg. 54]{SaffTotik} that $E$ is regular if and only if $U^{\mu_E}$ is continuous\footnote{Since $\supp(\mu_E)\subseteq\partial\Omega$ \cite[Cor. I.4.5]{SaffTotik}, it is again enough to check continuity of $U^{\mu_E}$ only on $\partial\Omega$.} in $\C$.

Often we use the concept of {\it balayage} of a measure (\cite[Sec. II.4]{SaffTotik}). Let $D$ be a domain (connected open set) with compact boundary $\partial D$ whose complement has positive capacity, and $\mu$ be a finite Borel measure with compact support in $D$. Then there exists a unique Borel measure $\widehat\mu$ supported on $\partial D$, with total mass is equal to that of $\mu$, whose potential $U^{\widehat\mu}$ is bounded on $\partial D$ and satisfies for some constant
$c(\mu;D)$
\begin{equation}
\label{eq:equalBal}
U^{\widehat\mu}(z) =    U^\mu(z)+c(\mu;D) \;\;\; \mbox{for q.e.} \;\;\; z\in\C\setminus D.
\end{equation}
Necessarily then, we have that $c(\mu;D)=0$ if $D$ is bounded and
$c(\mu;D)=\int g_D(t,\infty)d\mu(t)$ otherwise, where $g_D(\cdot,\infty)$ is the {\it Green function} for $D$ with pole at infinity. Equality in (\ref{eq:equalBal}) holds for all $z\in\C\setminus\overline D$ and also at all regular points of $\partial D$. The measure $\widehat\mu$ is called the balayage of $\mu$ onto $\partial D$.It has the property that
\begin{equation}
\label{balayageh}
\int h\,d\mu=\int h\,{d\widehat\mu}
\end{equation}
for any function $h$ which is harmonic in $D$ and continuous in $\overline{D}$ (including at infinity if $D$ is unbounded). From its defining properties $\widehat\mu$ has finite energy, therefore it cannot charge polar sets.

In analogy to the logarithmic case, one can define the {\it Green potential} of a positive measure $\mu$ supported in a domain $D$ with compact non-polar boundary. The only difference is now that, in (\ref{eq1sL}), the
logarithmic kernel $\log(1/|z-t|)$ gets replaced by  $g_D(z,t)$, the {\it Green function} for $D$ with pole at $t\in D$.  The Green potential relative to the domain $D$ of a finite positive Borel measure $\mu$ compactly supported in $D$ is given by
\[U_D^\mu(z)=\int g_D(z,t)\,d\mu(t).\]
It can be re-expressed in terms of the logarithmic potentials of $\mu$ and of its balayage $\widehat\mu$ onto $\partial D$ by the following formula \cite[Thm. II.4.7 and Thm. II.5.1]{SaffTotik}:
\begin{equation}
\label{eq:toRemind}
U^{\widehat\mu-\mu}(z) = c(\mu;D) - U_D^\mu(z), \;\;\; z\in D,
\end{equation}
where $c(\mu;D)$ was defined after equation (\ref{eq:equalBal}). Moreover, (\ref{eq:toRemind}) continues to hold at every regular point of $\partial D$; in particular, it holds q.e. on $\partial D$.

Exactly as in the logarithmic case, if $E$ is a compact nonpolar subset of $D$, there exists a unique measure $\mu_{(E,\partial D)}\in\Lm(E)$ that minimizes the Green energy among all measures in $\Lm(E)$. This measure is called the {\it Green equilibrium distribution} on $E$ relative to $D$. In addition, the Green equilibrium distribution satisfies
\begin{equation}
\label{eq:GreenEqual}
U_D^{\mu_{(E,\partial D)}}(z)=\frac{1}{\cp(E,\partial D)}, \;\;\; \mbox{for q.e.} \;\;\; z\in E,
\end{equation}
where $\cp(E,\partial D)$ is {\it Green (condenser) capacity} of $E$ relative to $D$ which is the reciprocal of the minimal Green energy among all measures in $\Lm(E)$. Moreover, equality in (\ref{eq:GreenEqual}) holds at all regular points of $E$.

\small\centering
\bibliographystyle{plain}
\bibliography{weakpade}

\begin{thebibliography}{10}

\bibitem{AAK71}
V.~M. Adamyan, D.~Z. Arov, and M.~G. Krein.
\newblock Analytic properties of {S}chmidt pairs for a {H}ankel operator on the
  generalized {S}chur-{T}akagi problem.
\newblock {\em Math. USSR Sb.}, 15:31--73, 1971.

\bibitem{Ank83}
A.~Ancona.
\newblock D\'emonstration d'une conjecture sur la capacit\'e et l'effilement.
\newblock {\em C. R. Acad. Sci. Paris S}, 297(7):393--395, 1983.

\bibitem{Ank_CTP84}
A.~Ancona.
\newblock Sur une conjecture concernant la capacit\'e et l'effilement.
\newblock In G.~Mokobodzi and D.~Pinchon, editors, {\em Colloque du Th\'eorie
  du Potentiel (Orsay, 1983)}, volume 1096 of {\em Lecture Notes in
  Mathematics}, pages 34--68, Springer-Verlag, Berlin, 1984.

\bibitem{Ap02}
A.~I. Aptekarev.
\newblock Sharp constant for rational approximation of analytic functions.
\newblock {\em Mat. Sb.}, 193(1):1--72, 2002.
\newblock English transl. in {\it {M}ath. {S}b.} 193(1-2):1--72, 2002.

\bibitem{AVA04}
A.~I. Aptekarev and W.~V. Assche.
\newblock Scalar and matrix {R}iemann-{H}ilbert approach to the strong
  asymptotics of {P}ad\'e approximants and complex orthogonal polynomials with
  varying weight.
\newblock {\em J. Approx. Theory}, 129:129--166, 2004.

\bibitem{BKT05}
L.~Baratchart, R.~K\"ustner, and V.~Totik.
\newblock Zero distribution via orthogonality.
\newblock {\em Ann. Inst. Fourier}, 55(5):1455--1499, 2005.

\bibitem{BMSW06}
L.~Baratchart, F.~Mandr\`ea, E.~B. Saff, and F.~Wielonsky.
\newblock 2-{D} inverse problems for the {L}aplacian: a meromorphic
  approximation approach.
\newblock {\em J. Math. Pures Appl.}, 86:1--41, 2006.

\bibitem{BS02}
L.~Baratchart and F.~Seyfert.
\newblock An ${L}^p$ analog of {AAK} theory for $p\geq2$.
\newblock {\em J. Funct. Anal.}, 191(1):52--122, 2002.

\bibitem{uBY1}
L.~Baratchart and M.~Yattselev.
\newblock Meromorphic approximants to complex {C}auchy transforms with polar
  singularities.
\newblock {\it Submitted for publication}.

\bibitem{BarrLLTorr95}
D.~Barrios, G.~L\'opez Lagomasino, and E.~Torrano.
\newblock Location of the zeros and asymptotics of polynomials satisfying three
  term recurrence relations with complex coefficients.
\newblock {\em Sb. Math.}, 186(5):629--659, 1995.

\bibitem{Bax61}
G.~Baxter.
\newblock A convergence equivalence related to polynomials orthogonal on the
  unit circle.
\newblock {\em Trans, Amer. Math. Soc.}, 79:471--487, 1961.

\bibitem{Deift}
P.~Deift.
\newblock {\em Orthogonal Polynomials and Random Matrices: a Riemann-Hilbert
  Approach}, volume~3 of {\em Courant Lectures in Mathematics}.
\newblock Amer. Math. Soc., Providence, RI, 2000.

\bibitem{DeVoreLorentz}
R.~A. DeVore and G.~G. Lorentz.
\newblock {\em Constructive Approximation}, volume 303 of {\em Grundlehren der
  Math. Wissenschaften}.
\newblock Springer-Verlag, Berlin, 1993.

\bibitem{Gon75b}
A.~A. Gonchar.
\newblock On the convergence of generalized {P}ad\'e approximants of
  meromorphic functions.
\newblock {\em Mat. Sb.}, 98(140):564--577, 1975.
\newblock English transl. in {\it {M}ath. {USSR} {S}b.} 27:503--514, 1975.

\bibitem{Gon75a}
A.~A. Gonchar.
\newblock On the convergence of {P}ad\'e approximants for some classes of
  meromorphic functions.
\newblock {\em Mat. Sb.}, 97(139):607--629, 1975.
\newblock English transl. in {\it {M}ath. {USSR} {S}b.} 26(4):555--575, 1975.

\bibitem{GL78}
A.~A. Gonchar and G.~L\'opez Lagomasino.
\newblock On {M}arkov's theorem for multipoint {P}ad\'e approximants.
\newblock {\em Mat. Sb.}, 105(4):512--524, 1978.
\newblock English transl. in {\it{M}ath. {USSR} {S}b.} 34(4):449--459, 1978.

\bibitem{GRakh87}
A.~A. Gonchar and E.~A. Rakhmanov.
\newblock Equilibrium distributions and the degree of rational approximation of
  analytic functions.
\newblock {\em Mat. Sb.}, 134(176)(3):306--352, 1987.
\newblock English transl. in {\it {M}ath. {USSR} Sbornik} 62(2):305--348, 1989.

\bibitem{thKus}
R.~K\"ustner.
\newblock {\em Asymptotic zero distribution of orthogonal polynomials with
  respect to complex measures having argument of bounded variation}.
\newblock PhD thesis, {U}niversity of {N}ice {S}ophia {A}ntipolis, Sophia
  Antipolis, France, 2003.
\newblock http://www.inria.fr/rrrt/tu-0784.html.

\bibitem{LL78}
G.~L\'opez Lagomasino.
\newblock Conditions for convergence of multipoint {P}ad\'e approximants for
  functions of {S}tieltjes type.
\newblock {\em Mat. Sb.}, 107(149):69--83, 1978.
\newblock English transl. in {\it{M}ath. {USSR} {S}b.} 35:363--376, 1979.

\bibitem{LL88}
G.~L\'opez Lagomasino.
\newblock Convergence of {P}ad\'e approximants of {S}tieltjes type meromorphic
  functions and comparative asymptotics for orthogonal polynomials.
\newblock {\em Mat. Sb.}, 136(178):206--226, 1988.
\newblock English transl. in {\it{M}ath. {USSR} {S}b.} 64:207--227, 1989.

\bibitem{LLMarVA95}
G.~L\'opez Lagomasino, F.~Marcell\'an, and W.~Van Assche.
\newblock Relative asymptotics for polynomials orthogonal with respect to a
  discrete {S}obolev inner product.
\newblock {\em Constr. Approx.}, 11(1):107--137, 1995.

\bibitem{Mag87}
A.~Magnus.
\newblock Toeplitz matrix techniques and convergence of complex weight {P}ad\'e
  approximants.
\newblock {\em J. Comput. Appl. Math.}, 19:23--38, 1987.

\bibitem{Mar95}
A.~A. Markov.
\newblock Deux d\'emonstrations de la convergence de certaines fractions
  continues.
\newblock {\em Acta Math.}, 19:93--104, 1895.

\bibitem{Nut90}
J.~Nuttall.
\newblock Pad\'e polynomial asymptotic from a singular integral equation.
\newblock {\em Constr. Approx.}, 6(2):157--166, 1990.

\bibitem{NS80}
J.~Nuttall and S.~R. Singh.
\newblock Orthogonal polynomials and {P}ad\'e approximants associated with a
  system of arcs.
\newblock {\em J. Approx. Theory}, 21:1--42, 1980.

\bibitem{Par86}
O.~G. Parfenov.
\newblock Estimates of the singular numbers of a {C}arleson operator.
\newblock {\em Mat. Sb.}, 131(171):501--518, 1986.
\newblock English. transl. in {\it {M}ath. {USSR} {S}b. } 59:497--514, 1988.

\bibitem{Rakh77b}
E.~A. Rakhmanov.
\newblock Convergence of diagonal {P}ad\'e approximants.
\newblock {\em Mat. Sb.}, 104(146):271--291, 1977.
\newblock English transl. in {\it{M}ath. {USSR} {S}b.} 33:243--260, 1977.

\bibitem{Rakh77a}
E.~A. Rakhmanov.
\newblock On the asymptotics of the ratio of orthogonal polynomials.
\newblock {\em Mat. Sb.}, 103(145):237--252, 1977.
\newblock English transl. in {\it{M}ath. {USSR} {S}b.} 32:199--213, 1977.

\bibitem{Ransford}
T.~Ransford.
\newblock {\em Potential Theory in the Complex Plane}, volume~28 of {\em London
  Mathematical Society Student Texts}.
\newblock Cambridge University Press, Cambridge, 1995.

\bibitem{SaffTotik}
E.~B. Saff and V.~Totik.
\newblock {\em Logarithmic Potentials with External Fields}, volume 316 of {\em
  Grundlehren der Math. Wissenschaften}.
\newblock Springer-Verlag, Berlin, 1997.

\bibitem{St89}
H.~Stahl.
\newblock On the convergence of generalized {P}ad\'e approximants.
\newblock {\em Constr. Approx.}, 5(2):221--240, 1989.

\bibitem{St97}
H.~Stahl.
\newblock The convergence of {P}ad\'e approximants to functions with branch
  points.
\newblock {\em J. Approx. Theory}, 91:139--204, 1997.

\bibitem{StahlTotik}
H.~Stahl and V.~Totik.
\newblock {\em General Orthogonal Polynomials}, volume~43 of {\em Encycl.
  Math.}
\newblock Cambridge University Press, Cambridge, 1992.

\bibitem{Suet00}
S.~P. Suetin.
\newblock Uniform convergence of {P}ad\'e diagonal approximants for
  hyperelliptic functions.
\newblock {\em Mat. Sb.}, 191(9):81--114, 2000.
\newblock English transl. in {\it {M}ath. {S}b.} 191(9):1339--1373, 2000.

\bibitem{Suet02}
S.~P. Suetin.
\newblock Approximation properties of the poles of diagonal {P}ad\'e
  approximants for certain generalizations of {M}arkov functions.
\newblock {\em Mat. Sb.}, 193(12):105--133, 2002.
\newblock English transl. in {\it {M}ath. {S}b.} 193(12):1837--1866, 2002.

\bibitem{Totik}
V.~Totik.
\newblock {\em Weighted Approximation with Varying Weights}, volume 1300 of
  {\em Lecture {N}otes in {M}ath.}
\newblock Springer-Verlag, Berlin, 1994.

\end{thebibliography}

\end{document}